\documentclass[12pt,oneside,reqno]{amsart}
\usepackage{mathrsfs}
\usepackage{graphics}
\usepackage{graphicx,color}
\usepackage{amssymb}
\usepackage{subfigure}
\usepackage{enumerate}
\newcommand{\dif}{\mathrm{d}}

\newcommand{\be}{\begin{eqnarray}}
\newcommand{\ee}{\end{eqnarray}}
\newcommand{\ce}{\begin{eqnarray*}}
\newcommand{\de}{\end{eqnarray*}}
\newtheorem{theorem}{Theorem}[section]
\newtheorem{lemma}[theorem]{Lemma}
\newtheorem{remark}[theorem]{Remark}
\newtheorem{definition}[theorem]{Definition}
\newtheorem{proposition}[theorem]{Proposition}
\newtheorem{Example}[theorem]{Example}
\newtheorem{corollary}[theorem]{Corollary}
\def\e{\varepsilon}
\def\t{\theta}

\def\[{{\Big[}}
\def\]{{\Big]}}
\def\<{{\langle}}
\def\>{{\rangle}}
\def\({{\Big(}}
\def\){{\Big)}}

\def\no{\nonumber}
\def\bt{\begin{theorem}}
\def\et{\end{theorem}}
\def\bl{\begin{lemma}}
\def\el{\end{lemma}}
\def\br{\begin{remark}}
\def\er{\end{remark}}
\def\bx{\begin{Example}}
\def\ex{\end{Example}}
\def\bd{\begin{definition}}
\def\ed{\end{definition}}
\def\bp{\begin{proposition}}
\def\ep{\end{proposition}}
\def\bc{\begin{corollary}}
\def\ec{\end{corollary}}

\def\cB{{\mathcal B}}
\def\cC{{\mathcal C}}

\def\cJ{{\mathcal J}}
\def\cK{{\mathcal K}}

\def\cM{{\mathcal M}}
\def\cN{{\mathcal N}}

\def\cR{{\mathcal R}}

\def\cW{{\mathcal W}}

\def\mE{{\mathbb E}}

\def\mH{{\mathbb H}}

\def\mM{{\mathbb M}}
\def\mN{{\mathbb N}}

\def\mP{{\mathbb P}}

\def\mR{{\mathbb R}}

\def\mX{{\mathbb X}}

\def\geq{\geqslant}
\def\leq{\leqslant}

\def\sF{{\mathscr F}}

\def\sX{{\mathscr X}}

\title{Effective Filtering for Multiscale Stochastic Dynamical Systems in Hilbert Spaces*}

\author{Huijie Qiao}

\dedicatory{School of Mathematics \& Jiangsu Provincial Key Laboratory of Networked Collective Intelligence,\\
Southeast University, Nanjing, Jiangsu 211189,  China\\
hjqiaogean@seu.edu.cn}
\thanks{{\it AMS Subject Classification(2010):} 60H15, 37D10, 70K70}

\thanks{{\it Keywords:} Multiscale systems in Hilbert spaces, random invariant manifolds, nonlinear filtering, dimensional reduction.}

\thanks{*This work was supported by NSF of China (No. 11001051, 11371352) and the Jiangsu Provincial Key Laboratory of Networked Collective Intelligence under Grant No. BM2017002.}

\subjclass{}

\date{}

\begin{document}

\allowdisplaybreaks

\begin{abstract}
In the paper, effective filtering for a type of slow-fast data assimilation systems in Hilbert spaces is considered. Firstly, the system is reduced to a system on a random invariant manifold. Secondly, nonlinear filtering of the origin system can be approximated by that of the reduction system. Finally, we apply the obtained result to an example.
\end{abstract}

\maketitle \rm

\section{Introduction}

Give a probability space $(\Omega,\sF,\mP)$ and two separable Hilbert spaces $\mH^1, \mH^2$ with the inner products $\<\cdot, \cdot\>_{\mH^1}, \<\cdot, \cdot\>_{\mH^2}$ and the norms $\|\cdot\|_{\mH^1}, \|\cdot\|_{\mH^2}$, respectively. Consider a stochastic slow-fast system on $\mH^1\times\mH^2$
\be\left\{\begin{array}{l}
\dot{x}^\e=A{x}^\e+F(x^\e,y^\e)+\sigma_1\dot{W_1}, \\
\dot{y}^\e=\frac{1}{\e}B{y}^\e+\frac{1}{\e}G(x^\e, y^\e)+\frac{\sigma_2}{\sqrt{\e}}\dot{W_2},
\label{slfasy}
\end{array}
\right.
\ee
where $A, B$ are two linear operators on $\mH^1, \mH^2$, respectively, and the interaction functions $F:\mH^1\times\mH^2\rightarrow\mH^1$ and $G:\mH^1\times\mH^2\rightarrow\mH^2$ are Borel measurable. Moreover,  $W_1, W_2$ are two-sided $\mH^1, \mH^2$-valued Brownian motions with covariance operators $K_1, K_2$ such that $trK_1<\infty, trK_2<\infty$, respectively, and mutually independent.  $\sigma_1$ and $\sigma_2$ are nonzero real noise intensities, and $\e$ is a small positive parameter representing the ratio of two time scales. The type of systems (\ref{slfasy}) have appeared in many fields, such as engineering and science(\cite{ZS}). For example, the climate evolution consists of fast atmospheric and slow oceanic dynamics, and state dynamic in electric power systems include fast- and slowly-varying elements. 

The research for systems (\ref{slfasy}) is various. Let us mention some referrences. Schmalfu{\ss}-Schneider \cite{Sch} observed the invariant manifold for systems (\ref{slfasy}) in finite dimensional Hilbert spaces $\mH^1, \mH^2$. When $\mH^1, \mH^2$ are infinite dimensional, and only the fast part contains a finite dimensional noise, Fu-Liu-Duan \cite{Fu} studied the invariant manifold of systems (\ref{slfasy}). Stochastic average of systems (\ref{slfasy}) is considered in \cite{RG, andrew}. By comparing the slow manifold reduction with the averageing principle reduction, we point out that the conditions of the former  are shorter and the reduction system of the former is easier to simulate on the computer. 

Fix a separable Hilbert space $\mH^3$ with the inner product $\<\cdot, \cdot\>_{\mH^3}$ and the norm $\|\cdot\|_{\mH^3}$. The nonlinear filtering  problem for the slow component $x_t^{\e}$ with respect to a $\mH^3$-valued observation process $\{r^{\e}_s, 0\leq s\leq t\}$ (See Subsection \ref{nonfil} in details)  is to evaluate the `filter' $\mE[\phi(x^{\e}_t)|\mathcal{R}^{\e}_t]$, where $\phi$ is a Borel measurable function such that $\mE|\phi(x^{\e}_t)|<\infty$ for $t\in[0,T]$, and $\mathcal{R}^{\e}_t$ is the $\sigma$-algebra generated by $\{r^{\e}_s, 0\leq s\leq t\}$. In essence, the nonlinear filtering  is to achieve the best estimate for the origin system
state, given only noisy observations for the system. It provides an algorithm for estimating a
state of a random dynamical system based on noisy measurements. Stochastic filtering  is important in many practical applications,  from inertial guidance of aircrafts
and spacecrafts to weather and climate prediction. 

When $\mH^1, \mH^2, \mH^3$ are finite dimensional, the nonlinear filtering problems of multi-scale systems have been widely studied. Let us recall some results. In \cite{Imkeller}, Imkeller-Namachchivaya-Perkowski-Yeong showed that the filter  $\mE[\phi(x^{\e}_t)|\mathcal{R}^{\e}_t]$ converges to the homogenized filter by double backward stochastic differential equations and asymptotic techniques. Recently, Papanicolaou-Spiliopoulos \cite{ps1} also studied this convergence problem by independent version technique  and then applied it to statistical inference. When jumps processes are added in the system (\ref{slfasy}), the author proved the convergence by weak convergence technique in \cite{q1}. Besides, in \cite{q2} and \cite{zqd}, the author and two coauthors reduced the system (\ref{slfasy}) to a system on a random invariant manifold, and showed that $\mE[\phi(x^{\e}_t)|\mathcal{R}^{\e}_t]$ converges to the filter of the reduction system. Thus, a new method to study the nonlinear filtering  problem for multiscale systems is offered.

For a general system $x^{\e}$ on a Hilbert space, that is to say, there is no fast component $y^{\e}$, its nonlinear filtering  problem has been studied by Sritharan \cite{sr} and Hobbs-Sritharan \cite{hs}. However, for  nonlinear filtering  problems of multiscale systems on Hilbert spaces, nowadays there are no related results. Moreover, the type of problems have appeared in applications.(cf. \cite{yz})

In the paper, we consider a nonlinear filtering for the system (\ref{slfasy}) in general Hilbert spaces by following up the line in \cite{q2} and \cite{zqd}. Firstly, the system is reduced to a system on a random invariant manifold. Moreover, our result covers the known result in \cite{ccl, Fu}. Secondly, the nonlinear filtering of the origin system can be approximated by that of the reduction system. And this result generalizes the result in \cite{q2}.

It is worthwhile to mention our conditions and techniques. Firstly, the linear operators $A$ and $B$ may be unbounded, which can contain usual differential operators. Secondly, a general method of constructing random invariant manifolds for fast-slow systems is to transform those origin systems to non-random systems by stationary solutions of auxiliary equations and then construct random invariant manifolds for the obtained systems (c.f.\cite{Sch, q2, zqd}, etc.). However,  the existence of stationary solutions usually need strong conditions. Here we directly construct a random invariant manifold in terms of the system (\ref{slfasy}) itself. Therefore, our conditions are weaker than that in \cite{Sch, q2, zqd}. Finally, since these stochastic evolution equations on random slow manifolds  have no Markov property, some techniques, such as the Zakai equations in \cite{Park1, Park2, Park3} and backward stochastic differential equations in \cite{Imkeller}, do not work. Therefore, we use an exponential martingale technique to treat these nonlinear filtering problems. 

\bigskip

This paper is arranged as follows. In Section \ref{pre},  we introduce basic concepts about random dynamical systems and random invariant manifolds. We present the framework of the probability space and the system in Section \ref{fra}. Our method for reduced filtering is placed in Section \ref{redu}. In Section \ref{filter}, the nonlinear filtering problem is introduced and the approximation theorem of the filtering is proved. And then, we apply the obtained result to an example in Section \ref{exam}. In Section \ref{Con}, we summarize all the results in the paper. Finally, the proof of Theorem \ref{solexiuni} is placed in Section \ref{App}.

The following convention will be used throughout the paper: $C$ with or without indices will denote
different positive constants (depending on the indices) whose values may change from one place to
another. 

\section{Preliminaries}\label{pre}

In the section, we introduce some notations and basic concepts in random dynamical systems.

\subsection{Notation and terminology}

Let $\mathscr{B}(\mH^1)$ be the Borel $\sigma$-algebra on $\mH^1$, and $\cB(\mH^1)$ be the set of all real-valued bounded Borel-measurable functions on $\mH^1$.  Let $\cC(\mH^1)$ be the set
of all real-valued continuous functions on $\mH^1$, and $\cC_b^1(\mH^1)$ denote the collection of all functions of $\cC(\mH^1)$ which are bounded and Lipschitz continuous. And set
\ce
\|\phi\|:=\max\limits_{x\in\mH^1}|\phi(x)|+\max\limits_{x_1\neq x_2}\frac{|\phi(x_1)-\phi(x_2)|}{\|x_1- x_2\|_{\mH^1}}, \quad \phi\in\cC_b^1(\mH^1).
\de

\subsection{Random dynamical systems (\cite{la})}

Let $(\Omega,\sF,\mP)$ be a probability space, and $(\theta_t)_{t\in\mR}$ a family of measurable
transformations from $\Omega$ to $\Omega$ satisfying for $s,t\in\mR$,
\be
\theta_0=1_\Omega, \quad \theta_{t+s}=\theta_t\circ\theta_s.
\label{sepr}
\ee
If for each $t\in\mR$,
$\theta_t$ preserves the probability measure $\mP$, i.e.,
$$
\theta_t^*\mP=\mP,
$$
$(\Omega,\sF,\mP; (\theta_t)_{t\in\mR})$ is called a metric dynamical system. 

\bd\label{rds} Let $(\mX,\sX)$ be a measurable space. A mapping
\ce
\varphi: \mR\times\Omega\times\mX\mapsto\mX, \quad
(t,\omega,y)\mapsto\varphi(t,\omega,y)
\de
is called a measurable random dynamical system (RDS), or in short,  a
cocycle, if these following properties hold:

(i) Measurability: $\varphi$ is
$\mathscr{B}(\mR)\otimes\sF\otimes\sX/\sX$-measurable,

(ii) Cocycle property: $\varphi(t,\omega)$ satisfies the following conditions
$$
\varphi(0,\omega)=id_{\mX},
$$
and for $\omega\in\Omega$ and all $s,t\in\mR$ 
$$
\varphi(t+s,\omega)=\varphi(t,\theta_s\omega)\circ\varphi(s,\omega),
$$

(iii) Continuity: $\varphi(t,\omega)$ is continuous for $t\in\mR$.
\ed

\subsection{Random invariant manifolds (\cite{Sch})}

Let $\varphi$ be a random dynamical system on the normed space $(\mX, \|\cdot\|_{\mX})$. 

A family of nonempty  sets $\cM=\{\cM(\omega)\}_{\omega\in\Omega}$ is called a {\it random 
set} in $\mX$ if for $\omega\in\Omega$, $\cM(\omega)$ is a closed set in $\mX$ and for every $y\in\mX$, the mapping
$$
\Omega\ni\omega\mapsto dist(y,\cM(\omega)):=\inf\limits_{x\in \cM(\omega)}\|x-y\|_{\mX}
$$
is measurable. Moreover, if $\cM$ satisfies 
$$
\varphi(t, \omega, \cM(\omega))\subset \cM(\t_t\omega), \quad t\geq 0, \quad \omega\in\Omega,
$$
$\cM$ is called (positively) invariant with respect to $\varphi$.

In the sequel, we consider a random set defined by a Lipschitz continuous function. Concretely speaking, we will define a mapping by $$
\Omega\times\mH^1 \ni (\omega,x)\mapsto H(\omega, x)\in\mH^2,
$$
where for all $\omega\in\Omega$, $H(\omega, x)$ is globally Lipschitzian in $x$  and for any
$x\in\mH^1$, $\omega\mapsto H(\omega, x)$ is a $\mH^2$-valued random variable.  Thus, set
$$
\cM(\omega) :=\{(x, H(\omega, x))|x\in\mH^1\},
$$
and then $\cM$ is a random set in $\mH^1\times\mH^2$ (\cite[Lemma 2.1]{Sch}). Moreover, the random set
$\cM(\omega)$ is called a {\it Lipschitz random invariant manifold} if it is (positively) invariant with respect to some random dynamical system.

\section{Framework}\label{fra}

In the section, we present some results which will be applied in the following sections.

Let $\Omega^1:=C_0(\mR, \mH^1)$ be the collection of all strongly continuous functions $f: \mR\rightarrow\mH^1$ with $f(0)=0$. And then it is equipped with the compact-open topology. Let $\sF^1$ be its Borel
$\sigma$-algebra and $\mP^1$ the distribution of $W_1$ on $\Omega^1$.  Set
$$
\t_t^1\omega_1(\cdot):=\omega_1(\cdot+t)-\omega_1(t), \quad \omega_1\in\Omega^1, \quad t\in\mR,
$$
and then $(\theta^1_t)_{t\in\mR}$ satisfy (\ref{sepr}). Moreover, by the property of $\mP^1$ we obtain that $(\Omega^1, \sF^1, \mP^1, \t_t^1)$ is a metric dynamical system.  Next, set $\Omega^2:=C_0(\mR, \mH^2)$. And then we define $\sF^2, \mP^2, \t_t^2$ by the similar means to $\sF^1, \mP^1, \t_t^1$. Thus, $(\Omega^2, \sF^2, \mP^2, \t_t^2)$ becomes another metric dynamical system.  Set
$$
\Omega:=\Omega^1\times\Omega^2, ~\sF:=\sF^1\times\sF^2, ~\mP:=\mP^1\times\mP^2,~
\t_t:=\t_t^1\times\t_t^2,
$$
and then $(\Omega, \sF, \mP, \t_t)$ is a metric dynamical system that is used in the sequel.

\bigskip

Consider the slow-fast system (\ref{slfasy}) on $\mH^1\times\mH^2$, i.e.
\ce\left\{\begin{array}{l}
\dot{x}^\e=A{x}^\e+F(x^\e,y^\e)+\sigma_1 \dot{W_1}, \\
\dot{y}^\e=\frac{1}{\e}B{y}^\e+\frac{1}{\e}G(x^\e, y^\e)+\frac{\sigma_2}{\sqrt{\e}}\dot{W_2}.
\end{array}
\right.
\de
We make the following hypotheses:

\medspace

\begin{enumerate}[($\bf{H_1}$)] 
\item There exists a $\gamma_1\geq0$ such that
\be
\|e^{At}\|\leq e^{-\gamma_1t}, \quad t\leq0,
\label{opecon1}
\ee
where $\|e^{At}\|$ stands for the norm of the operator $e^{At}$, and $\{e^{At}, t\geq0\}$ is a strongly continuous semigroup on $\mH^1$ and 
\be
\|e^{At}\|\leq 1,  \quad t\geq0.
\label{opecon2}
\ee
\end{enumerate}
\begin{enumerate}[($\bf{H_2}$)] 
\item There exists a $\gamma_2>0$ such that for any $y\in\mH^2$,
$$
\<By,y\>_{\mH^2}\leq -\gamma_2\|y\|^2_{\mH^2}.
$$
\end{enumerate}
\begin{enumerate}[($\bf{H_3}$)] 
\item There exists a positive constant $L$ such that for all $(x_1, y_1), (x_2, y_2)\in \mH^1\times\mH^2$
\ce
\|F(x_1, y_1)-F(x_2, y_2)\|_{\mH^1}\leq L(\|x_1-x_2\|_{\mH^1}+\|y_1-y_2\|_{\mH^2}),
\de
and
\ce
\|G(x_1, y_1)-G(x_2, y_2)\|_{\mH^2}\leq L(\|x_1-x_2\|_{\mH^1}+\|y_1-y_2\|_{\mH^2}),
\de
and $F(0,0)=G(0,0)=0$.
\end{enumerate}
\begin{enumerate}[($\bf{H_4}$)] 
\item
$$
\gamma_2>L.
$$
\end{enumerate}

\br
By ($\bf{H_2}$), we know that $\frac{B}{\e}$ generates a strongly continuous semigroup $\{e^{\frac{B}{\e}t}, t\geq0\}$ on $\mH^2$ and 
\be
\|e^{\frac{B}{\e}t}\|\leq e^{-\frac{\gamma_2}{\e}t}, \quad t\geq0.
\label{deccon}
\ee
And ($\bf{H_3}$) admits us to obtain that for any $x\in\mH^1, y\in\mH^2$ 
\be
\|F(x, y)\|_{\mH^1}\leq L(\|x\|_{\mH^1}+\|y\|_{\mH^2}), \quad \|G(x, y)\|_{\mH^2}\leq L(\|x\|_{\mH^1}+\|y\|_{\mH^2}).
\label{lingro}
\ee
\er

In the following, we give the definition of mild solutions to the system (\ref{slfasy}) and then state that the system (\ref{slfasy}) has a unique mild solution which generates a RDS. Let $\mH:=\mH^1\times\mH^2$ with the norm $\|z\|_{\mH}=\|x\|_{\mH^1}+\|y\|_{\mH^2}$ for $z=(x,y)\in\mH$. Let $\cC([a,b], \mH)$ be the collection of strongly continuous functions on $[a,b]$ with values in $\mH$.

\bd\label{mildef}
Let $s\in\mR, T>0$ and $z_0=(x_0, y_0)\in\mH$. $z^\e(t)\equiv z^\e(t, s, \omega; z_0)$ is said to be a mild solution to the system (\ref{slfasy}) on the interval $(s, s+T]$ if (i) $z^\e(s)=z_0$, (ii) $z^\e(t)$ satisfies the following integral equation
\ce
z^\e(t)=\left(\begin{array}{cccc}
x_t^\e\\
y_t^\e
\end{array}\right)=\left(\begin{array}{cccc}e^{A(t-s)}x_0+\int_s^te^{A(t-r)}F(x_r^\e,y_r^\e)dr+\int_s^te^{A(t-r)}\sigma_1\dif W_1(r)\\
e^{\frac{B}{\e}(t-s)}y_0+\int_s^te^{\frac{B}{\e}(t-r)}\frac{1}{\e}G(x_r^\e,y_r^\e)dr+\int_s^te^{\frac{B}{\e}(t-r)}\frac{\sigma_2}{\sqrt{\e}}\dif W_2(r)\end{array}\right)
\de
for $t\in[s,s+T]$ and $\omega\in\Omega$, (iii) $z^\e(t)$ belongs to $\cC([s,s+T], \mH)$.
\ed

\bt\label{solexiuni}
Suppose that ($\bf{H_1}$)-($\bf{H_4}$) are satisfied. Let $s\in\mR, T>0$ and $z_0=(x_0, y_0)\in\mH$. Then the system (\ref{slfasy}) has a unique mild solution $z^\e(t, s, \omega; z_0)$ for $t\in[s,s+T]$ and $\omega\in\Omega$. Moreover, set $\varphi^\e(t,\omega)z_0:=z^\e(t, 0, \omega; z_0), t\geq 0$, and then $\varphi^\e(t,\omega)$ is a RDS.
\et

The above theorem is similar to Theorem 3.3 in \cite{ccl}. To the readers' convenience, we give a proof in the Appendix.

\section{Reduction on random invariant manifolds}\label{redu}

In the section, we prove that the origin system can be approximated by the reduction system on a random invariant manifold.

\subsection{Random invariant manifolds}

In the subsection, we prove that the system (\ref{slfasy}) has a random invariant manifold. Firstly, we  introduce some spaces. Let
\ce
&&\cC^{1,-}_{\frac{\mu}{\e},s}:=\left\{\phi\in\cC((-\infty,s], \mH^1): \sup\limits_{t\leq s}e^{\frac{\mu}{\e} (t-s)}\|\phi(t)\|_{\mH^1}<\infty\right\},\\
&&\cC^{2,-}_{\frac{\mu}{\e},s}:=\left\{\phi\in\cC((-\infty,s], \mH^2): \sup\limits_{t\leq s}e^{\frac{\mu}{\e} (t-s)}\|\phi(t)\|_{\mH^2}<\infty\right\},\\
&&\cC^{1}_{\frac{\mu}{\e}}:=\left\{\phi\in\cC(\mR, \mH^1): \sup\limits_{t\in\mR}e^{\frac{\mu}{\e} t}\|\phi(t)\|_{\mH^1}<\infty\right\},\\
&&\cC^{2}_{\frac{\mu}{\e}}:=\left\{\phi\in\cC(\mR, \mH^2): \sup\limits_{t\in\mR}e^{\frac{\mu}{\e} t}\|\phi(t)\|_{\mH^2}<\infty\right\},
\de
where $\mu$ is a positive constant and satisfies $\gamma_2-L >\mu>0$. Let $\cC^-_{\frac{\mu}{\e},s}:=\cC^{1,-}_{\frac{\mu}{\e},s}\times\cC^{2,-}_{\frac{\mu}{\e},s}$ with the norm $\|z\|_{\cC^-_{\frac{\mu}{\e},s}}=\sup\limits_{t\leq s}e^{\frac{\mu}{\e}(t-s)}\|z(t)\|_{\mH}$ for $z\in\cC^-_{\frac{\mu}{\e},s}$, and $\cC_{\frac{\mu}{\e}}:=\cC^{1}_{\frac{\mu}{\e}}\times\cC^{2}_{\frac{\mu}{\e}}$ with the norm $\|z\|_{\cC_{\frac{\mu}{\e}}}=\sup\limits_{t\in\mR}e^{\frac{\mu}{\e}t}\|z(t)\|_{\mH}$ for $z\in\cC_{\frac{\mu}{\e}}$.

\bl\label{equi}
Suppose that $(\bf{H_1})$--$(\bf{H_4})$ are satisfied. Let $s\in\mR$ and $x_0\in\mH^1$. Then there exists a $\e_0>0$ such that for $0<\e\leq\e_0$, the following integral equation has a unique solution $\bar{z}_t^\e=(\bar{x}_t^\e,\bar{y}_t^\e)$ for $t\leq s$:
\be
\bar{z}_t^\e=\left(\begin{array}{cccc}
\bar{x}_t^\e\\
\bar{y}_t^\e
\end{array}\right)=\left(\begin{array}{cccc}e^{A(t-s)}x_0-\int_t^s e^{A(t-r)}F(\bar{x}_r^\e,\bar{y}_r^\e)dr-\int_t^s e^{A(t-r)}\sigma_1\dif W_1(r)\\
\int_{-\infty}^te^{\frac{B}{\e}(t-r)}\frac{1}{\e}G(\bar{x}_r^\e,\bar{y}_r^\e)dr+\int_{-\infty}^te^{\frac{B}{\e}(t-r)}\frac{\sigma_2}{\sqrt{\e}}\dif W_2(r)\end{array}\right), t\leq s,
\label{ineq}
\ee
\ce
\bar{x}_s^\e=x_0.
\de
\el
\begin{proof}
Set for $\bar{z}^\e=(\bar{x}^\e,\bar{y}^\e)\in\cC^-_{\frac{\mu}{\e},s}$
\be
&&\cK(\bar{z}^\e)(t):=\left(\begin{array}{cccc}
\cK_1(\bar{z}^\e)(t)\\
\cK_2(\bar{z}^\e)(t)
\end{array}\right):=\left(\begin{array}{cccc}e^{A(t-s)}x_0-\int_t^se^{A(t-r)}F(\bar{x}_r^\e,\bar{y}_r^\e)dr-\int_t^se^{A(t-r)}\sigma_1\dif W_1(r)\\
\int_{-\infty}^te^{\frac{B}{\e}(t-r)}\frac{1}{\e}G(\bar{x}_r^\e,\bar{y}_r^\e)dr+\int_{-\infty}^te^{\frac{B}{\e}(t-r)}\frac{\sigma_2}{\sqrt{\e}}\dif W_2(r)\end{array}\right), \no\\
&&\qquad\qquad\qquad\qquad\qquad\qquad\qquad\qquad t\leq s,
\label{kop}
\ee
and then $\cK$ is well defined on $\cC^-_{\frac{\mu}{\e},s}$. Indeed, we calculate that for $\bar{z}^\e=(\bar{x}^\e,\bar{y}^\e)\in\cC^-_{\frac{\mu}{\e},s}$,
\ce
\sup\limits_{t\leq s}e^{\frac{\mu}{\e} (t-s)}\|e^{A(t-s)}x_0\|_{\mH^1}\leq\sup\limits_{t\leq s}e^{\frac{\mu}{\e} (t-s)}e^{-\gamma_1 (t-s)}\|x_0\|_{\mH^1}\leq\|x_0\|_{\mH^1},
\de
\be
\sup\limits_{t\leq s}e^{\frac{\mu}{\e} (t-s)}\left\|\int_t^se^{A(t-r)}F(\bar{x}_r^\e,\bar{y}_r^\e)dr\right\|_{\mH^1}&\leq&\sup\limits_{t\leq s}e^{\frac{\mu}{\e} (t-s)}\int_t^s e^{-\gamma_1 (t-r)}L(\|\bar{x}_r^\e\|_{\mH^1}+\|\bar{y}_r^\e\|_{\mH^2})dr\no\\
&\leq&L\left(\sup\limits_{t\leq s}e^{\frac{\mu}{\e} (t-s)}\|\bar{z}_t^\e\|_{\mH}\right)\sup\limits_{t\leq s}\int_t^s e^{(\frac{\mu}{\e}-\gamma_1) (t-r)}dr\no\\
&=&\frac{\e L}{\mu-\e\gamma_1}\left(\sup\limits_{t\leq s}e^{\frac{\mu}{\e} (t-s)}\|\bar{z}_t^\e\|_{\mH}\right), \quad a.s.,
\label{kf}
\ee
and
\be
\sup\limits_{t\leq s}e^{\frac{\mu}{\e} (t-s)}\left\|\int_{-\infty}^t e^{\frac{B}{\e}(t-r)}\frac{1}{\e}G(\bar{x}_r^\e,\bar{y}_r^\e)dr\right\|_{\mH^2}
&\leq&\frac{1}{\e}\sup\limits_{t\leq s}e^{\frac{\mu}{\e} (t-s)}\int_{-\infty}^te^{-\frac{\gamma_2}{\e}(t-r)}L(\|\bar{x}_r^\e\|_{\mH^1}+\|\bar{y}_r^\e\|_{\mH^2})dr\no\\
&\leq&\frac{L}{\e}\left(\sup\limits_{t\leq s}e^{\frac{\mu}{\e} (t-s)}\|\bar{z}_t^\e\|_{\mH}\right)\int_{-\infty}^te^{(\frac{\mu}{\e}-\frac{\gamma_2}{\e})(t-r)}dr\no\\
&=&\frac{L}{\gamma_2-\mu}\left(\sup\limits_{t\leq s}e^{\frac{\mu}{\e} (t-s)}\|\bar{z}_t^\e\|_{\mH}\right), \quad a.s..
\label{kg}
\ee
By \cite[Proposition 3.1]{cs}, it holds that 
\ce
&&\sup\limits_{t\leq s}e^{\frac{\mu}{\e} (t-s)}\left\|\int_t^s e^{A(t-r)}\sigma_1\dif W_1(r)\right\|_{\mH^1}<\infty,\quad a.s.,\\
&&\sup\limits_{t\leq s}e^{\frac{\mu}{\e} (t-s)}\left\|\int_{-\infty}^te^{\frac{B}{\e}(t-r)}\frac{\sigma_2}{\sqrt{\e}}\dif W_2(r)\right\|_{\mH^2}<\infty, \quad a.s..
\de

In the following, we show that $\cK$ is contractive. For $\bar{z}^{\e,1}, \bar{z}^{\e,2}\in\cC^-_{\frac{\mu}{\e},s}$, by the same deduction to (\ref{kf}) (\ref{kg}), one can obtain that 
\ce
&&\sup\limits_{t\leq s}e^{\frac{\mu}{\e}(t-s)}\left\|\cK_1(\bar{z}^{\e,1})(t)-\cK_1(\bar{z}^{\e,2})(t)\right\|_{\mH^1}\leq\frac{\e L}{\mu-\e\gamma_1}\left(\sup\limits_{t\leq s}e^{\frac{\mu}{\e} (t-s)}\|\bar{z}_t^{\e,1}-\bar{z}_t^{\e,2}\|_{\mH}\right), \\
&&\sup\limits_{t\leq s}e^{\frac{\mu}{\e} (t-s)}\left\|\cK_2(\bar{z}^{\e,1})(t)-\cK_2(\bar{z}^{\e,2})(t)\right\|_{\mH^2}\leq\frac{L}{\gamma_2-\mu}\left(\sup\limits_{t\leq s}e^{\frac{\mu}{\e} (t-s)}\|\bar{z}_t^{\e,1}-\bar{z}_t^{\e,2}\|_{\mH}\right).
\de
Thus, we get that
\ce
\sup\limits_{t\leq s}e^{\frac{\mu}{\e} (t-s)}\left\|\cK(\bar{z}^{\e,1})(t)-\cK(\bar{z}^{\e,2})(t)\right\|_{\mH}&\leq&\sup\limits_{t\leq s}e^{\frac{\mu}{\e} (t-s)}\left\|\cK_1(\bar{z}^{\e,1})(t)-\cK_1(\bar{z}^{\e,2})(t)\right\|_{\mH^1}\\
&&+\sup\limits_{t\leq s}e^{\frac{\mu}{\e} (t-s)}\left\|\cK_2(\bar{z}^{\e,1})(t)-\cK_2(\bar{z}^{\e,2})(t)\right\|_{\mH^2}\\
&\leq&\left(\frac{\e L}{\mu-\e\gamma_1}+\frac{L}{\gamma_2-\mu}\right)\left(\sup\limits_{t\leq s}e^{\frac{\mu}{\e} (t-s)}\|\bar{z}_t^{\e,1}-\bar{z}_t^{\e,2}\|_{\mH}\right).
\de
Therefore, in order to prove that $\cK$ is contractive, one only needs to show that 
$$
\frac{\e L}{\mu-\e\gamma_1}+\frac{L}{\gamma_2-\mu}<1.
$$
Note that $\gamma_2-L>\mu>0$, and then 
$$
\frac{L}{\gamma_2-\mu}<1.
$$
Thus, there exists a $\e_0>0$ such that for any $0<\e\leq\e_0$,
$$
\frac{\e L}{\mu-\e\gamma_1}+\frac{L}{\gamma_2-\mu}<1,
$$
and furthermore $\cK$ is contractive. So, Eq.(\ref{ineq}) has a unique solution. The proof is completed.
\end{proof}

Next, we study a property of the solution to Eq.(\ref{ineq}). First of all, for $u\in(-\infty, s)$, Eq.(\ref{ineq}) can be rewritten as 
\ce
\left(\begin{array}{cccc}
\bar{x}_{t}^\e\\
\bar{y}_{t}^\e
\end{array}\right)=\left(\begin{array}{cccc}e^{A(t-u)}\bar{x}_{u}^\e+\int_u^te^{A(t-r)}F(\bar{x}_{r}^\e,\bar{y}_{r}^\e)dr+\int_u^te^{A(t-r)}\sigma_1\dif W_1(r)\\
e^{\frac{B}{\e}(t-u)}\bar{y}_{u}^\e+\int_u^te^{\frac{B}{\e}(t-r)}\frac{1}{\e}G(\bar{x}_{r}^\e,\bar{y}_{r}^\e)dr+\int_u^te^{\frac{B}{\e}(t-r)}\frac{\sigma_2}{\sqrt{\e}}\dif W_2(r)\end{array}\right), u\leq t\leq s.
\de
That is, the dynamic of Eq.(\ref{ineq}) is the same as that of the system (\ref{slfasy}). Thus, by Theorem \ref{solexiuni} it holds that for $r\leq s$
\be
\bar{x}_{r}^\e(\t_u\omega)=\bar{x}_{r+u}^\e(\omega), \quad \bar{y}_{r}^\e(\t_u\omega)=\bar{y}_{r+u}^\e(\omega), \quad u\leq s-r.
\label{apro}
\ee

\br\label{invmancon}
Set
\ce
H^{\e,s}(\omega,x_0):=\bar{y}_s^\e=\int_{-\infty}^se^{\frac{B}{\e}(s-r)}\frac{1}{\e}G(\bar{x}_{r}^\e,\bar{y}_{r}^\e)dr+\int_{-\infty}
^se^{\frac{B}{\e}(s-r)}\frac{\sigma_2}{\sqrt{\e}}\dif W_2(r),
\de
and then for $x_0^1, 
x_0^2\in\mH^1$, it holds that 
\be
\|H^{\e,s}(\omega,x^1_0)-H^{\e,s}(\omega,x^2_0)\|_{\mH^2}\leq\frac{Le^{\frac{\mu}{\e} s}}{(\gamma_2-\mu)\left[1-\(\frac{\e L}{\mu-\e\gamma_1}+\frac{L}{\gamma_2-\mu}\)\right]}\|x_0^1-x_0^2\|_{\mH^1},
\label{lipcon}
\ee
where we use (\ref{kf}) (\ref{kg}). Besides, it follows from (\ref{apro}) that for $t\leq s$,
\be
H^{\e,s}(\t_{t-s}\omega,x_0)&=&\int_{-\infty}^se^{\frac{B}{\e}(s-r)}\frac{1}{\e}G\left(\bar{x}_{r}^\e(\t_{t-s}\omega),\bar{y}_{r}^\e(\t_{t-s}\omega)\right)dr\no\\
&&+\int_{-\infty}
^se^{\frac{B}{\e}(s-r)}\frac{\sigma_2}{\sqrt{\e}}\dif W_2(r)(\t_{t-s}\omega)\no\\
&=&\int_{-\infty}^te^{\frac{B}{\e}(t-r)}\frac{1}{\e}G(\bar{x}_{r}^\e,\bar{y}_{r}^\e)dr+\int_{-\infty}
^te^{\frac{B}{\e}(t-r)}\frac{\sigma_2}{\sqrt{\e}}\dif W_2(r)\no\\
&=&\bar{y}_{t}^\e.
\label{stapro}
\ee
\er

\bt\label{inma} (A random invariant manifold) \\
Assume that $(\bf{H_1})$--$(\bf{H_4})$ are satisfied. Let $z_0=(x_0, y_0)\in\mH$. Then for $0<\e\leq\e_0$, $\varphi^\e$ has a random invariant manifold
$$
\cM^\e(\omega)=\left\{\big(x, H^{\e,0}(\omega,x)\big), x\in\mH^1\right\},
$$
where for $\omega\in\Omega$, the Lipschitz constant of $H^{\e,0}(\omega,x)$ is bounded by
$$
\frac{L}{(\gamma_2-\mu)\left[1-\(\frac{\e L}{\mu-\e\gamma_1}+\frac{L}{\gamma_2-\mu}\)\right]}.
$$
Moreover, $\cM^\e$ is exponentially attracting in the following sense: for any $z_0=(x_0, y_0)$, there exists a $\tilde{z}_0=(\tilde{x}_0,\tilde{y}_0)\in\cM^\e(\omega)$ such that
\be
\|\varphi^\e(t,\omega)z_0-\varphi^\e(t,\omega)\tilde{z}_0\|_{\mH}\leq \frac{e^{-\frac{\mu}{\e} t}}{1-M}\((2+2M)\|z_0\|_{\mH}+2 R(\omega)\), ~t\geq0,
\label{esti}
\ee
where 
\ce
M:= \frac{\e L}{\mu-\e\gamma_1}+\frac{L}{\gamma_2-\mu},
\de
\ce
R(\omega):=\sup\limits_{t\leq 0}e^{\frac{\mu}{\e} t}\left\|\int_t^0 e^{A(t-r)}\sigma_1\dif W_1(r)\right\|_{\mH^1}+\sup\limits_{t\leq 0}e^{\frac{\mu}{\e} t}\left\|\int_{-\infty}^te^{\frac{B}{\e}(t-r)}\frac{\sigma_2}{\sqrt{\e}}\dif W_2(r)\right\|_{\mH^2}.
\de 
\et
\begin{proof}
First of all, set 
\ce
\cM^\e(\omega):=\left\{\big(x, H^{\e,0}(\omega,x)\big), x\in\mH^1\right\},
\de
and then by (\ref{lipcon}), it holds that $\cM^\e(\omega)$ is a random set. And based on (\ref{stapro}) one can justify that $\cM^\e(\omega)$ is invariant with respect to $\varphi^\e$. Therefore, $\cM^\e(\omega)$ is a Lipschitz random invariant manifold with respect to $\varphi^\e$. 

Next, we prove (\ref{esti}). For this, set 
\ce
z_t^\e=\left(\begin{array}{cccc}
x_t^\e\\
y_t^\e
\end{array}\right):=\left\{\begin{array}{cccc}
\(x_0, (I-|t|B)^{-1}y_0\), \qquad t\leq0,\\
\varphi^\e(t,\omega)(x_0,y_0), \qquad\quad\quad\quad t>0,
\end{array}\right.
\de
and
\ce
Z_0(t):=\left\{\begin{array}{cccc}
-z_t^\e+\cK(z^\e)(t), \qquad\qquad\qquad\qquad\qquad\quad t\leq0,\\
-(e^{At}x_0,e^{\frac{B}{\e}t}y_0)+(e^{At}x_0,e^{\frac{B}{\e}t}\cK_2(z^\e)(0)), \quad t>0,
\end{array}\right.
\de
where $\cK$ is defined in (\ref{kop}). And then we consider the following integral equation
\be
\left(\begin{array}{cccc}
X_t^\e\\
Y_t^\e
\end{array}\right)
=Z_0(t)+\left(\begin{array}{cccc}-\int_t^\infty e^{A(t-r)}\[F(x_r^\e+X_r^\e, y_r^\e+Y_r^\e)-F(x_r^\e, y_r^\e)\]dr\\
\int_{-\infty}^t e^{\frac{B}{\e}(t-r)}\frac{1}{\e}\[G(x_r^\e+X_r^\e, y_r^\e+Y_r^\e)-G(x_r^\e, y_r^\e)\]dr\end{array}\right), t\in\mR.
\label{ineq3}
\ee
For $Z^\e:=(X^\e, Y^\e)\in\cC_{\frac{\mu}{\e}}$, set
\ce
\cR(Z^\e)(t)&&:=
\left(\begin{array}{cccc}
\cR_1(Z^\e)(t)\\
\cR_2(Z^\e)(t)
\end{array}\right)\\
&&:=Z_0(t)+\left(\begin{array}{cccc}-\int_t^\infty e^{A(t-r)}\[F(x_r^\e+X_r^\e, y_r^\e+Y_r^\e)-F(x_r^\e, y_r^\e)\]dr\\
\int_{-\infty}^t e^{\frac{B}{\e}(t-r)}\frac{1}{\e}\[G(x_r^\e+X_r^\e, y_r^\e+Y_r^\e)-G(x_r^\e, y_r^\e)\]dr\end{array}\right),
\de
and then $\cR: \cC_{\frac{\mu}{\e}}\mapsto\cC_{\frac{\mu}{\e}}$ is well defined. Indeed, for $Z^\e=(X^\e, Y^\e)\in\cC_{\frac{\mu}{\e}}$, by $(\bf{H_1})$--$(\bf{H_4})$ we compute
\be
\sup\limits_{t\in\mR}e^{\frac{\mu}{\e} t}\|Z_0(t)\|_{\mH}&\leq&\sup\limits_{t\leq0}e^{\frac{\mu}{\e} t}\|Z_0(t)\|_{\mH}+\sup\limits_{t>0}e^{\frac{\mu}{\e} t}\|Z_0(t)\|_{\mH}\no\\
&\leq&\sup\limits_{t\leq0}e^{\frac{\mu}{\e} t}(\|z_t^\e\|_{\mH}+\|\cK(z^\e)(t)\|_{\mH})+\sup\limits_{t>0}e^{\frac{\mu}{\e} t}e^{-\frac{\gamma_2}{\e}t}\(\|y_0\|_{\mH^2}+\|\cK_2(z^\e)(0)\|_{\mH^2}\)\no\\
&\leq&\sup\limits_{t\leq0}e^{\frac{\mu}{\e} t}\|z_t^\e\|_{\mH}+\|x_0\|_{\mH^1}+M\sup\limits_{t\leq0}e^{\frac{\mu}{\e} t}\|z_t^\e\|_{\mH}+R(\omega)+\|y_0\|_{\mH^2}\no\\
&&+\|\cK_2(z^\e)(0)\|_{\mH^2}\no\\
&\leq&\sup\limits_{t\leq0}e^{\frac{\mu}{\e} t}\|z_t^\e\|_{\mH}+\|x_0\|_{\mH^1}+M\sup\limits_{t\leq0}e^{\frac{\mu}{\e} t}\|z_t^\e\|_{\mH}+R(\omega)+\|y_0\|_{\mH^2}\no\\
&&+M\sup\limits_{t\leq0}e^{\frac{\mu}{\e} t}\|z_t^\e\|_{\mH}+R(\omega)\no\\
&\leq&(1+2M)\sup\limits_{t\leq0}e^{\frac{\mu}{\e} t}\|z_t^\e\|_{\mH}+\|x_0\|_{\mH^1}+\|y_0\|_{\mH^2}+2 R(\omega)\no\\
&\leq&(2+2M)(\|x_0\|_{\mH^1}+\|y_0\|_{\mH^2})+2 R(\omega)\no\\
&=&(2+2M)\|z_0\|_{\mH}+2 R(\omega).
\label{esti0}
\ee
Besides, it follows from $(\bf{H_1})$--$(\bf{H_4})$ that
\be
&&\sup\limits_{t\in\mR}e^{\frac{\mu}{\e} t}\int_t^\infty e^{-\gamma_1(t-r)}\|F(x_r^\e+X_r^\e, y_r^\e+Y_r^\e)-F(x_r^\e, y_r^\e)\|_{\mH^1}dr\no\\
&\leq&L\(\sup\limits_{t\in\mR}e^{\frac{\mu}{\e} t}\|Z_t^\e\|_{\mH}\)\sup\limits_{t\in\mR}\int_t^\infty e^{(\frac{\mu}{\e}-\gamma_1)(t-r)}dr\no\\
&\leq&\frac{\e L}{\mu-\e\gamma_1}\(\sup\limits_{t\in\mR}e^{\frac{\mu}{\e} t}\|Z_t^\e\|_{\mH}\),
\label{esti1}
\ee
and
\be
&&\frac{1}{\e}\sup\limits_{t\in\mR}e^{\frac{\mu}{\e} t}\int_{-\infty}^t e^{-\frac{\gamma_2}{\e}(t-r)}\|G(x_r^\e+X_r^\e, y_r^\e+Y_r^\e)-G(x_r^\e, y_r^\e)\|_{\mH^2}dr\no\\
&\leq&\frac{L}{\e}\(\sup\limits_{t\in\mR}e^{\frac{\mu}{\e} t}\|Z_t^\e\|_{\mH}\)\sup\limits_{t\in\mR}\int_{-\infty}^t e^{(\frac{\mu}{\e}-\frac{\gamma_2}{\e})(t-r)}dr\no\\
&\leq&\frac{L}{\gamma_2-\mu}\(\sup\limits_{t\in\mR}e^{\frac{\mu}{\e} t}\|Z_t^\e\|_{\mH}\).
\label{esti2}
\ee
Thus, by combining (\ref{esti0}) (\ref{esti1}) with (\ref{esti2}), one can get that
\ce
\sup\limits_{t\in\mR}e^{\frac{\mu}{\e} t}\|\cR(Z^\e)(t)\|_{\mH}\leq\sup\limits_{t\in\mR}e^{\frac{\mu}{\e} t}\|\cR_1(Z^\e)(t)\|_{\mH^1}+\sup\limits_{t\in\mR}e^{\frac{\mu}{\e} t}\|\cR_2(Z^\e)(t)\|_{\mH^2}<\infty.
\de

Next, for $Z^{\e,1}, Z^{\e,2}\in\cC_{\frac{\mu}{\e}}$, by the similar deduction to (\ref{esti1}) (\ref{esti2}) we know that
\ce
\sup\limits_{t\in\mR}e^{\frac{\mu}{\e} t}\|\cR_1(Z^{\e,1})(t)-\cR_1(Z^{\e,2})(t)\|_{\mH^1}&\leq&\frac{\e L}{\mu-\e\gamma_1}\(\sup\limits_{t\in\mR}e^{\frac{\mu}{\e} t}\|Z_t^{\e,1}-Z_t^{\e,2}\|_{\mH}\),\\
\sup\limits_{t\in\mR}e^{\frac{\mu}{\e} t}\|\cR_2(Z^{\e,1})(t)-\cR_2(Z^{\e,2})(t)\|_{\mH^2}&\leq&\frac{L}{\gamma_2-\mu}\(\sup\limits_{t\in\mR}e^{\frac{\mu}{\e} t}\|Z_t^{\e,1}-Z_t^{\e,2}\|_{\mH}\).
\de
Thus, one can have that
\ce
 \sup\limits_{t\in\mR}e^{\frac{\mu}{\e} t}\|\cR(Z^{\e,1})(t)-\cR(Z^{\e,2})(t)\|_{\mH}&\leq&\sup\limits_{t\in\mR}e^{\frac{\mu}{\e} t}\|\cR_1(Z^{\e,1})(t)-\cR_1(Z^{\e,2})(t)\|_{\mH^1}\\
 &&+\sup\limits_{t\in\mR}e^{\frac{\mu}{\e} t}\|\cR_2(Z^{\e,1})(t)-\cR_2(Z^{\e,2})(t)\|_{\mH^2}\\
 &\leq&\(\frac{\e L}{\mu-\e\gamma_1}+\frac{L}{\gamma_2-\mu}\)\(\sup\limits_{t\in\mR}e^{\frac{\mu}{\e} t}\|Z_t^{\e,1}-Z_t^{\e,2}\|_{\mH}\).
\de
So, for $0<\e\leq\e_0$, $\cR: \cC_{\frac{\mu}{\e}}\mapsto \cC_{\frac{\mu}{\e}}$ is contractive. That is, Eq.(\ref{ineq3}) has a unique solution denoted as $Z^\e=(X^\e, Y^\e)$. Moreover, 
\ce
\sup\limits_{t\in\mR}e^{\frac{\mu}{\e} t}\|Z_t^\e\|_{\mH}\leq \frac{1}{1-M}\((2+2M)\|z_0\|_{\mH}+2 R(\omega)\)
\de
and then 
\be
\|Z_t^\e\|_{\mH}\leq \frac{e^{-\frac{\mu}{\e} t}}{1-M}\((2+2M)\|z_0\|_{\mH}+2 R(\omega)\)
, \quad t\geq 0.
\label{esti3}
\ee

Set
\ce
\tilde{z}_t^\e:=\left(\begin{array}{cccc}
\tilde{x}_t^\e\\
\tilde{y}_t^\e\end{array}\right):=\left(\begin{array}{cccc}
x_t^\e\\
y_t^\e\end{array}\right)+\left(\begin{array}{cccc}
X_t^\e\\
Y_t^\e\end{array}\right),
\de
and then by simple calculation, it holds that $\tilde{z}_t^\e=(\tilde{x}_t^\e,\tilde{y}_t^\e)$ solves uniquely the following equation
\ce
\tilde{z}_t^\e=\left\{\begin{array}{cccc}\cK(\tilde{z}^\e)(t),\qquad\quad t\leq0,\\
\varphi^\e(t,\omega)(\tilde{x}_0^\e, \tilde{y}_0^\e), \quad t>0.
\end{array}\right.
\de
Thus, by Lemma \ref{equi}, we know that $\tilde{z}_t^\e$ solves Eq.(\ref{ineq}) for $t\leq0$. In particular, $\tilde{y}_0^\e=H^{\e,0}(\omega, \tilde{x}_0^\e)$, which yields that $(\tilde{x}_0^\e, \tilde{y}_0^\e)\in\cM^\e(\omega)$. So, we take $\tilde{z}_0=(\tilde{x}_0^\e, \tilde{y}_0^\e)$. Since $\tilde{z}_t^\e=\varphi^\e(t,\omega)(\tilde{x}_0^\e, \tilde{y}_0^\e)$ for $t>0$, $\tilde{z}_t^\e=\varphi^\e(t,\omega)\tilde{z}_0$ for $t>0$.  Note that $\tilde{x}_t^\e-x_t^\e=X_t^\e, \tilde{y}_t^\e-y_t^\e=Y_t^\e$. Thus, by (\ref{esti3}), it holds that
\ce
\|\varphi^\e(t,\omega)z_0-\varphi^\e(t,\omega)\tilde{z}_0\|_{\mH}&=&\|(x_t^\e, y_t^\e)-(\tilde{x}_t^\e, \tilde{y}_t^\e)\|_{\mH}=\|Z_t^\e\|_{\mH}\\
&\leq&\frac{e^{-\frac{\mu}{\e} t}}{1-M}\((2+2M)\|z_0\|_{\mH}+2 R(\omega)\), \quad t\geq 0.
\de
The proof is completed.
\end{proof}

\subsection{A reduced system on the random invariant manifold $\cM^\e$}

In the subsection, we prove that there exists a reduced system on the random invariant manifold $\cM^\e$ such that it will approximate the original system (\ref{slfasy}) for sufficiently long time. 

By Theorem \ref{inma}, we can obtain the following result.

\bt\label{Reduction-Theo} (A reduced system on the random invariant manifold) \\
Assume that $(\bf{H_1})$--$(\bf{H_4})$ hold. Let $z_0=(x_0, y_0)\in\mH$. Then for $0<\e\leq\e_0$ and the system (\ref{slfasy}) with the initial value $z^\e(0)=z_0$, there exists the following system on the random invariant manifold $\cM^\e$:
\be\left\{\begin{array}{l}
\dot{\tilde{x}}^\e=A\tilde{x}^\e+F\left(\tilde{x}^\e,\tilde{y}^\e\right)+\sigma_1 \dot{W_1},\\
\tilde{y}^\e=H^{\e,0}( \theta_{\cdot}\omega,
\tilde{x}^\e),
\end{array}
\right.
\label{redsys}
\ee
such that for almost all $\omega$,
\be
\|z^\e(t, \omega)-\tilde{z}^\e(t, \omega)\|_{\mH}&\leq&\frac{e^{-\frac{\mu}{\e} t}}{1-M}\((2+2M)\|z_0\|_{\mH}+2 R(\omega)\), ~ t\geq 0,
\label{orrees}
\ee
where $\tilde{z}^\e(t)=(\tilde{x}^\e(t), \tilde{y}^\e(t))$ is the solution of the system (\ref{redsys}) with
the initial value $\tilde{z}^\e(0)=(\tilde{x}^\e_0, H^{\e,0}(\omega, \tilde{x}^\e_0))$.
\end{theorem}

\br
By (\ref{orrees}), we know that when $\e$ is enough small or $t$ is sufficiently large, the system (\ref{redsys}) will approximate the system \eqref{slfasy}.
\er

\section{An approximate filter on the invariant manifold}\label{filter}

In the section we introduce nonlinear filtering problems for the system (\ref{slfasy}) and the reduced
system (\ref{redsys}) on the random  invariant manifold, and then study their relation.

\subsection{Nonlinear filtering problems}\label{nonfil}

In the subsection we introduce nonlinear filtering problems for the system (\ref{slfasy}) and the reduced system (\ref{redsys}).

Let $\{\beta_i(t,\omega)\}_{i\geq1}$ be a family of mutually independent one-dimensional Brownian motions on $(\Omega, \mathscr{F}, \mathbb{P})$. Construct a cylindrical Brownian motion on $\mH^3$ with respect to $(\Omega, \mathscr{F}, \mathbb{P})$ by
\begin{equation*}
W_3(t):=W_3(t, \omega):=\sum_{i=1}^\infty\beta_i(t,\omega)e_i, \quad \omega\in\Omega,\,\, t\in[0,\infty),
\end{equation*}
where $\{e_i\}_{i\geq1}$ is a complete orthonormal basis for $\mathbb{H}^3$.  It is easy to justify
that the covariance operator of the cylindrical Brownian motion $W_3$ is the identity operator $I$ on $\mathbb{H}^3$. Note that $W_3$ is not a process on $\mathbb{H}^3$. It is convenient to realize $W_3$ as a continuous process on an enlarged Hilbert space $\tilde{\mathbb{H}}^3$, the completion of $\mathbb{H}^3$ under the inner product 
$$
\<x,y\>_{\tilde{\mathbb{H}}^3}:=\sum\limits_{i=1}^\infty 2^{-i}\<x,e_i\>_{\mH^3}\<y,e_i\>_{\mH^3}, \quad x, y\in\mathbb{H}^3.
$$
Note that here $W_3$ may be either independent of $W_1$ and $W_2$, or dependent on $W_1$ and $W_2$( see\cite{ccc}). 

In the following, fix a Borel measurable function $h(x,y):\mH^1\times\mH^2\rightarrow\mH^3$. For $h$, we make the following hypothesis:
\begin{enumerate}[($\bf{H_5}$)] 
\item There exists a $C_h>0$ such that $\sup\limits_{(x,y)\in\mH^1\times\mH^2}\|h(x,y)\|_{\mH^3}\leq C_h$ and $h(x,y)$ is Lipschitz continuous in $(x,y)$ whose Lipschitz constant is denoted by $\|h\|_{Lip}$.
\end{enumerate}

Next, for $T>0$, an observation system  is given by
\ce
r^{\e}_t= W_3(t)+\int_0^t h(x^\e_s, y^\e_s) \dif s, \quad t\in[0,T].
\de
Under the assumption $(\bf{H_5})$, $r^{\e}$ is well defined.  Set
$$
(\Gamma_t^\e)^{-1}:=\exp\left\{-\int_0^t \<h(x_s^\e, y_s^\e),\dif W_3(s)\>_{\tilde{\mathbb{H}}^3}-\frac12\int_0^t\|h(x_s^\e, y_s^\e)\|_{\mH^3}^2\dif
s\right\},
$$
and then by \cite[Proposition 10.17]{dpz}, $(\Gamma_t^\e)^{-1}$ is an exponential martingale under $\mP$. Thus, \cite[Theorem 10.14]{dpz} admits us to obtain that $r^{\e}$ is a cylindrical Brownian motion under a new probability measure $\mP^\e$ via
$$
\frac{\dif \mP^\e}{\dif \mP}:=(\Gamma^\e_T)^{-1}.
$$

Rewrite $\Gamma_t^\e$ as
$$
\Gamma_t^\e=\exp\left\{\int_0^t \<h(x_s^\e, y_s^\e),\dif r^{\e}_s\>_{\tilde{\mathbb{H}}^3}-\frac12\int_0^t\|h(x_s^\e, y_s^\e)\|_{\mH^3}^2\dif
s\right\},
$$
and define
 $$
 \rho_t^\e(\phi) :=\mE^\e[\phi(x_t^\e)\Gamma^\e_t|\mathcal{R}_t^\e], \quad \phi\in \cB(\mH^1),
 $$
 where  $\mE^\e$ stands for the expectation under $\mP^\e$, $\mathcal{R}_t^\e \triangleq\sigma(r_s^\e:
 0\leq s \leq t) \vee \cN$ and $\cN$ is the collection of all $\mP$-measure zero sets. And set
 \ce
 \pi_t^\e(\phi) := \mE[\phi(x_t^\e)|\mathcal{R}_t^\e], \quad \phi\in \cB(\mH^1),
 \de
and then by the Kallianpur-Striebel formula it holds that
\ce
\pi^{\e}_t(\phi)=\frac{\rho^{\e}_t(\phi)}{\rho^{\e}_t(1)}.
\de
Here $ \rho_t^\e$ is
 called the nonnormalized filtering of $x_t^\e$ with respect to $\mathcal{R}_t^\e$, and $ \pi_t^\e$ is called the normalized filtering of $x_t^\e$ with respect to $\mathcal{R}_t^\e$, or the nonlinear filtering problem for $x_t^\e$ with respect to $\mathcal{R}_t^\e$.

Besides, we rewrite the reduced system (\ref{redsys}) as
\ce
\dot{\tilde{x}}^\e=A\tilde{x}^\e+\tilde{F}^\e(\omega, \tilde{x}^\e)+\sigma_1\dot{W_1},
\de
where $\tilde{F}^\e(\omega, x):=F(x, H^{\e,0}( \theta_{\cdot}\omega, x))$, and study the nonlinear filtering
problem for $\tilde{x}^\e$. Set
\ce
&&\tilde{h}^\e(\omega, x):=h(x, H^{\e,0}( \theta_{\cdot}\omega, x)),\\
&&\tilde{\Gamma}^\e_t:=\exp\left\{\int_0^t\<\tilde{h}^\e(\omega,\tilde{x}^\e_s),\dif
r^{\e}_s\>_{\tilde{\mathbb{H}}^3}-\frac12\int_0^t\|\tilde{h}^\e(\omega,\tilde{x}^\e_s)\|_{\mH^3}^2\dif s\right\},
\de
and then by \cite[Proposition 10.17]{dpz}, $\tilde{\Gamma}^\e_t$ is an exponential martingale under $\mP^\e$. Thus, we define 
\ce
&&\tilde{\rho}_t^\e(\phi) :=\mE^\e[\phi(\tilde{x}^\e_t)\tilde{\Gamma}^\e_t|\mathcal{R}_t^\e],\\
&&\tilde{\pi}_t^\e(\phi):=\frac{\tilde{\rho}_t^\e(\phi)}{\tilde{\rho}_t^\e(1)}, \quad \phi\in \cB(\mH^1),
\de
and then study the relation between $\pi^{\e}_t$ and $\tilde{\pi}^\e_t$.

\subsection{The relation between $\pi^{\e}_t$ and $\tilde{\pi}^\e_t$}

In the subsection we investigate the relation of $\pi^{\e}_t$ and $\tilde{\pi}^\e_t$ for $\e$ small enough. Let us start with two estimations.

\bl\label{es1}
Assume that $(\bf{H_1})$--$(\bf{H_5})$ are satisfied. Then 
$$
\mE^{\e}\left|\tilde{\rho}^{\e}_t(1)\right|^{-p}\leq\exp\left\{\left(\frac{p^2}{2}+\frac{p}{2}\right)C_h^2T\right\}, \quad t\in[0,T], \quad
p>0.
$$
\el
\begin{proof}
By the Jensen inequality it holds  that
$$
\mE^\e\left|\tilde{\rho}^{\e}_t(1)\right|^{-p}=\mE^\e\left|\mE^\e[\tilde{\Gamma}^\e_t|\mathcal{R}_t^\e]\right|^{-p}\leq\mE^\e\left[\mE^\e[|\tilde{\Gamma}^\e_t|^{-p}|\mathcal{R}_t^\e]\right]=\mE^\e[|\tilde{\Gamma}^\e_t|^{-p}].
$$
So, by  the definition of $\tilde{\Gamma}^\e_t$ we know that 
\ce
\mE^\e[|\tilde{\Gamma}^\e_t|^{-p}]&=&\mE^\e\left[\exp\left\{-p\int_0^t\<\tilde{h}^\e(\omega,\tilde{x}^\e_s), \dif
r^{\e}_s\>_{\tilde{\mathbb{H}}^3}+\frac{p}{2}\int_0^t\|\tilde{h}^\e(\omega,\tilde{x}^\e_s)\|_{\mH^3}^2\dif s\right\}\right]\\
&=&\mE^\e\Bigg[\exp\left\{-p\int_0^t\<\tilde{h}^\e(\omega,\tilde{x}^\e_s), \dif
r^{\e}_s\>_{\tilde{\mathbb{H}}^3}-\frac{p^2}{2}\int_0^t\|\tilde{h}^\e(\omega,\tilde{x}^\e_s)\|_{\mH^3}^2\dif s\right\}\\
&&\bullet\exp\left\{\left(\frac{p^2}{2}+\frac{p}{2}\right)\int_0^t\|\tilde{h}^\e(\omega,\tilde{x}^\e_s)\|_{\mH^3}^2\dif
s\right\}\Bigg]\\
&\leq&\mE^\e\Bigg[\exp\left\{-p\int_0^t\<\tilde{h}^\e(\omega,\tilde{x}^\e_s), \dif
r^{\e}_s\>_{\tilde{\mathbb{H}}^3}-\frac{p^2}{2}\int_0^t\|\tilde{h}^\e(\omega,\tilde{x}^\e_s)\|_{\mH^3}^2\dif s\right\}\Bigg]\\
&&\bullet\exp\left\{\left(\frac{p^2}{2}+\frac{p}{2}\right)C_h^2T\right\}\\
&=&\exp\left\{\left(\frac{p^2}{2}+\frac{p}{2}\right)C_h^2T\right\},
\de
where the last step is based on the fact that  $\exp\left\{-p\int_0^t\<\tilde{h}^\e(\omega,\tilde{x}^\e_s), \dif
r^{\e}_s\>_{\tilde{\mathbb{H}}^3}-\frac{p^2}{2}\int_0^t\|\tilde{h}^\e(\omega,\tilde{x}^\e_s)\|_{\mH^3}^2\dif s\right\}$ is an exponential martingale under $\mP^\e$. The proof is completed.
\end{proof}

\bl\label{es2}
Under $(\bf{H_1})$--$(\bf{H_5})$, it holds that  for $0<\e\leq\e_0$ and $\phi\in \cC^1_b(\mH^1)$,
\ce
\mE^{\e}\left|\rho^{\e}_t(\phi)-\tilde{\rho}^{\e}_t(\phi)\right|^p&\leq&\|\phi\|^{p}\frac{C}{\left[1-M\right]^p}\(\mE^\e\((2+2M)\|z_0\|_{\mH}+2 R(\omega)\)^{2p}\)^{1/2}\\
&&\qquad\qquad \times \(e^{-\frac{\mu}{\e} tp}+\frac{\e}{\mu p}\), \quad t\in[0,T],  \quad p>2,
\de
where $C>0$ is a constant independent of $\e$.
\el
\begin{proof}
For $\phi\in \cC^1_b(\mH^1)$ and $p>2$, by the definition of $\rho^{\e}_t(\phi), \tilde{\rho}^{\e}_t(\phi)$ it holds that
\be
\mE^\e\left|\rho^{\e}_t(\phi)-\tilde{\rho}^{\e}_t(\phi)\right|^{p}&=&\mE^\e\left|\mE^\e[\phi(x_t^\e)\Gamma^\e_t|\mathcal{R}_t^\e]-\mE^\e[\phi(\tilde{x}^\e_t)\tilde{\Gamma}^\e_t|\mathcal{R}_t^\e]\right|^{p}\no\\
&=&\mE^\e\left|\mE^\e[\phi(x_t^\e)\Gamma^\e_t-\phi(\tilde{x}^\e_t)\tilde{\Gamma}^\e_t|\mathcal{R}_t^\e]\right|^{p}\no\\
&\leq&\mE^\e\left[\mE^\e\left[\left|\phi(x_t^\e)\Gamma^\e_t-\phi(\tilde{x}^\e_t)\tilde{\Gamma }^\e_t\right|^{p}\bigg|\mathcal{R}_t^\e\right]\right]\no\\
&=&\mE^\e\left[\left|\phi(x_t^\e)\Gamma^\e_t-\phi(\tilde{x}^\e_t)\tilde{\Gamma}^\e_t\right|^{p}\right]\no\\
&\leq&2^{p-1}\mE^\e\left[\left|\phi(x_t^\e)\Gamma^\e_t-\phi(\tilde{x}^\e_t)\Gamma^\e_t\right|^{p}\right]\no\\
&&+2^{p-1}\mE^\e\left[\left|\phi(\tilde{x}_t^\e)\Gamma^\e_t-\phi(\tilde{x}^\e_t)\tilde{\Gamma}^\e_t\right|^{p}\right]\no\\
&=:&J_1+J_2.
\label{i1i2}
\ee

To $J_1$, by the H\"older inequality, we know that
\begin{equation}
\label{11}
\begin{array}{rcl}
J_1&\leq&\displaystyle
2^{p-1}(\mE^\e\left[\left|\phi(x_t^\e)-\phi(\tilde{x}^\e_t)\right|^{2p}\right])^{1/2}
(\mE^\e\left|\Gamma^\e_t\right|^{2p})^{1/2} \\[1ex]
&\leq&2^{p-1}\|\phi\|^{p}(\mE^\e\left\|x_t^\e-\tilde{x}^\e_t\right\|_{\mH^1}^{2p})^{1/2}\Bigg(\mE^\e\exp\left\{2p\int_0^t
\<h(x_s^\e, y_s^\e), \dif r^{\e}_s\>_{\tilde{\mathbb{H}}^3}-\frac{(2p)^2}{2}\int_0^t\|h(x_s^\e, y_s^\e)\|_{\mH^3}^2\dif s\right\}\\[2ex]
&& \bullet \exp\left\{\frac{(2p)^2}{2}\int_0^t\|h(x_s^\e, y_s^\e)\|_{\mH^3}^2\dif s-\frac{2p}{2}\int_0^t\|h(x_s^\e, y_s^\e)\|_{\mH^3}^2\dif s\right\}
\Bigg)^{1/2}\\[2ex]
&\leq&2^{p-1}\|\phi\|^{p}\frac{e^{-\frac{\mu}{\e} tp}}{\left[1-M\right]^p}\(\mE^\e\((2+2M)\|z_0\|_{\mH}+2 R(\omega)\)^{2p}\)^{1/2}e^{p(2p-1)C_h^2T/2},
\end{array}
\end{equation}
where the last step is based on Theorem \ref{Reduction-Theo} and the fact that the process\\
$\exp\left\{2p\int_0^t
\<h(x_s^\e, y_s^\e), \dif r^{\e}_s\>_{\tilde{\mathbb{H}}^3}-\frac{(2p)^2}{2}\int_0^t\|h(x_s^\e, y_s^\e)\|_{\mH^3}^2\dif s\right\}$ is an exponential martingale under $\mP^\e$.

Next, for $J_2$, it holds that
\ce
J_2\leq2^{p-1}\|\phi\|^{p}\mE^\e\left[\left|\Gamma^\e_t-\tilde{\Gamma}^\e_t\right|^{p}\right].
\de
Based on the It\^o formula, $\Gamma^\e_t$ and $\tilde{\Gamma}^\e_t$ solve the following equations, respectively,
\ce
\Gamma^\e_t=1+\int_0^t\Gamma^\e_s \<h(x_s^\e, y_s^\e), \dif r_s^\e\>_{\tilde{\mathbb{H}}^3}, \quad
\tilde{\Gamma}^\e_t=1+\int_0^t\tilde{\Gamma}^\e_s \<\tilde{h}^\e(\omega,\tilde{x}^\e_s), \dif r_s^\e\>_{\tilde{\mathbb{H}}^3}.
\de
So, it follows from the BDG inequality and the H\"older inequality that
\ce
\mE^\e\left[\left|\Gamma^\e_t-\tilde{\Gamma}^\e_t\right|^{p}\right]&=&\mE^\e\left[\left|\int_0^t\left<\Gamma^\e_s
h(x_s^\e, y_s^\e)-\tilde{\Gamma}^\e_s \tilde{h}^\e(\omega,\tilde{x}^\e_s),\dif
r_s^\e\right>_{\tilde{\mathbb{H}}^3}\right|^{p}\right]\\
&\leq&\mE^\e\left[\int_0^t\left\|\Gamma^\e_s h(x_s^\e, y_s^\e)-\tilde{\Gamma}^\e_s
\tilde{h}^\e(\omega,\tilde{x}^\e_s)\right\|_{\mH^3}^2\dif s\right]^{p/2}\\
&\leq&T^{p/2-1}\int_0^t\mE^\e\left\|\Gamma^\e_s h(x_s^\e, y_s^\e)-\tilde{\Gamma}^\e_s
\tilde{h}^\e(\omega,\tilde{x}^\e_s)\right\|_{\mH^3}^{p}\dif s\\
&\leq&2^{p-1}T^{p/2-1}\int_0^t\mE^\e\left\|\Gamma^\e_s h(x_s^\e, y_s^\e)-\Gamma^\e_s
\tilde{h}^\e(\omega,\tilde{x}^\e_s)\right\|_{\mH^3}^{p}\dif s\\
&&+2^{p-1}T^{p/2-1}\int_0^t\mE^\e\left\|\Gamma^\e_s
\tilde{h}^\e(\omega,\tilde{x}^\e_s)-\tilde{\Gamma}^\e_s
\tilde{h}^\e(\omega,\tilde{x}^\e_s)\right\|_{\mH^3}^{p}\dif s\\
&=:&J_{21}+J_{22}.
\de
For $J_{21}$, by the similar deduction to $J_1$ we have
\ce
J_{21}&\leq&2^{p-1}T^{p/2-1}\int_0^t\|h\|_{Lip}^{p}\frac{e^{-\frac{\mu}{\e} sp}}{\left[1-M\right]^p}\(\mE^\e\((2+2M)\|z_0\|_{\mH}+2 R(\omega)\)^{2p}\)^{1/2}\\
&&\qquad\qquad\qquad \times e^{p(2p-1)C_h^2T/2}\dif s\\
&=&2^{p-1}T^{p/2-1}\|h\|_{Lip}^{p}\frac{1}{\left[1-M\right]^p}\(\mE^\e\((2+2M)\|z_0\|_{\mH}+2 R(\omega)\)^{2p}\)^{1/2}\\
&&\qquad\qquad\qquad \times e^{p(2p-1)C_h^2T/2}\frac{\e}{\mu p}[1-e^{-\frac{\mu}{\e} t p}].
\de
And for $J_{22}$, it follows from the bounded property of $h$ that
\ce
J_{22}\leq2^{p-1}T^{p/2-1}C_h^p\int_0^t\mE^\e\left|\Gamma^\e_s-\tilde{\Gamma}^\e_s
\right|^{p}\dif s.
\de
So,
\ce
\mE^\e\left[\left|\Gamma^\e_t-\tilde{\Gamma}^\e_t\right|^{p}\right]&\leq&
C\cdot \frac{1}{\left[1-M\right]^p}\frac{\e}{\mu p}\(\mE^\e\((2+2M)\|z_0\|_{\mH}+2 R(\omega)\)^{2p}\)^{1/2}\\
&&+C\int_0^t\mE^\e\left|\Gamma^\e_s-\tilde{\Gamma}^\e_s
\right|^{p}\dif s,
\de
where the constant $C>0$ is independent of $\e, \mu$. By the Gronwall inequality it holds that
$$
\mE^\e\left[\left|\Gamma^\e_t-\tilde{\Gamma}^\e_t\right|^{p}\right]\leq
C\cdot \frac{1}{\left[1-M\right]^p}\frac{\e}{\mu p}\(\mE^\e\((2+2M)\|z_0\|_{\mH}+2 R(\omega)\)^{2p}\)^{1/2}.
$$
Thus,
\be
J_2\leq2^{p-1}\|\phi\|^{p}C\cdot \frac{1}{\left[1-M\right]^p}\frac{\e}{\mu p}\(\mE^\e\((2+2M)\|z_0\|_{\mH}+2 R(\omega)\)^{2p}\)^{1/2}.
\label{i2}
\ee

Finally, combining (\ref{i1i2}) with (\ref{11}) and (\ref{i2}), we obtain that
\ce
\mE^\e\left|\rho^{\e}_t(\phi)-\tilde{\rho}^{\e}_t(\phi)\right|^{p}&\leq&2^{p-1}\|\phi\|^{p}\frac{e^{-\frac{\mu}{\e} tp}}{\left[1-M\right]^p}\(\mE^\e\((2+2M)\|z_0\|_{\mH}+2 R(\omega)\)^{2p}\)^{1/2}\\
&&\qquad\qquad \times e^{p(2p-1)C_h^2T/2}\\
&&+2^{p-1}\|\phi\|^{p}C\cdot \frac{1}{\left[1-M\right]^p}\frac{\e}{\mu p}\\
&&\qquad\qquad \times\(\mE^\e\((2+2M)\|z_0\|_{\mH}+2 R(\omega)\)^{2p}\)^{1/2}\\
&\leq&\|\phi\|^{p}\frac{C}{\left[1-M\right]^p}\(\mE^\e\((2+2M)\|z_0\|_{\mH}+2 R(\omega)\)^{2p}\)^{1/2}\\
&&\qquad\qquad \times \(e^{-\frac{\mu}{\e} tp}+\frac{\e}{\mu p}\).
\de
This proves the lemma.
\end{proof}

Now, we are ready to state and prove the main result in the paper. First of all, we give out two concepts used in the proof of Theorem \ref{filcon}.

\bd\label{sepapo}
For the set $\mM\subset\cC_b^1(\mH^1)$, if the convergence $\lim\limits_{n\rightarrow\infty}\phi(x_n)=\phi(x), \forall \phi\in \mM$, for some $x, x_n\in\mH^1$, implies  that $\lim\limits_{n\rightarrow\infty} x_n =x$,  it is said that $\mM$ strongly separates points in $\mH^1$.
\ed

\bd\label{condet}
For the set $\mN\subset\cC_b^1(\mH^1)$, if $\mu_n$ and $\mu$ are probability measures on $\mathscr{B}(\mH^1)$, such that $\lim\limits_{n\rightarrow\infty} \int_{\mH^1}\phi\,\dif\mu_n = \int_{\mH^1}\phi\,\dif\mu$ for any $\phi\in \mN$, then $\mu_n$ converges weakly to $\mu$, it is said that $\mN$ is convergence determining for the topology of weak convergence of probability measures.
\ed

\bt\label{filcon} (Approximation by the reduced filter on the invariant manifold)\\
Under  $(\bf{H_1})$--$(\bf{H_5})$,  there exists a positive constant $C$ such that for $0<\e<\e_0$ and $\phi\in \cC_b^1(\mH^1)$
\ce
\mE|\pi^{\e}_t(\phi)- \tilde{\pi}_t^\e(\phi)|^p&\leq&\|\phi\|^{p}\frac{C}{\left[1-M\right]^p}\(\mE\((2+2M)\|z_0\|_{\mH}+2 R(\omega)\)^{16p}\)^{1/16}\\
&&\qquad\qquad \times \(e^{-\frac{4\mu}{\e} tp}+\frac{\e}{4\mu p}\)^{1/4}, \quad t\in[0,T],  \quad p>2.
\de
Thus, for the  distance $d(\cdot, \cdot)$  in the space of probability measures that induces the
weak convergence, the following approximation holds:
\ce
\mE [d(\pi^{\e}_t, \tilde{\pi}_t^\e)] &\leq&
\frac{C^{1/p}}{\left[1-M\right]}\(\mE\((2+2M)\|z_0\|_{\mH}+2 R(\omega)\)^{16p}\)^{1/16p}\\
&&\qquad\qquad \times \(e^{-\frac{4\mu}{\e} tp}+\frac{\e}{4\mu p}\)^{1/4p}.
\de
\et
\begin{proof}
For $\phi\in \cC^1_b(\mH^1)$, the H\"older inequality, Lemma \ref{es1} and Lemma \ref{es2} admit us to obtain that
\ce
\mE|\pi^{\e}_t(\phi)-
\tilde{\pi}_t^\e(\phi)|^{p}&=&\mE\left|\frac{\rho^{\e}_t(\phi)-\tilde{\rho}^{\e}_t(\phi)}{\tilde{\rho}^{\e}_t(1)}-\pi^{\e}_t(\phi)\frac{\rho^{\e}_t(1)-\tilde{\rho}^{\e}_t(1)}{\tilde{\rho}^{\e}_t(1)}\right|^{p}\\
&\leq&2^{p-1}\mE\left|\frac{\rho^{\e}_t(\phi)-\tilde{\rho}^{\e}_t(\phi)}{\tilde{\rho}^{\e}_t(1)}\right|^{p}+2^{p-1}\mE\left|\pi^{\e}_t(\phi)\frac{\rho^{\e}_t(1)-\tilde{\rho}^{\e}_t(1)}{\tilde{\rho}^{\e}_t(1)}\right|^{p}\\
&\leq&2^{p-1}\left(\mE\left|\rho^{\e}_t(\phi)-\tilde{\rho}^{\e}_t(\phi)\right|^{2p}\right)^{1/2}\left(\mE\left|\tilde{\rho}^{\e}_t(1)\right|^{-2p}\right)^{1/2}\\
&&+2^{p-1}\|\phi\|^{p}\left(\mE\left|\rho^{\e}_t(1)-\tilde{\rho}^{\e}_t(1)\right|^{2p}\right)^{1/2}\left(\mE\left|\tilde{\rho}^{\e}_t(1)\right|^{-2p}\right)^{1/2}\\
&=&2^{p-1}\left(\mE^{\e}\left|\rho^{\e}_t(\phi)-\tilde{\rho}^{\e}_t(\phi)\right|^{2p}\Gamma_T^{\e}\right)^{1/2}\left(\mE^{\e}\left|\tilde{\rho}^{\e}_t(1)\right|^{-2p}\Gamma_T^{\e}\right)^{1/2}\\
&&+2^{p-1}\|\phi\|^{p}\left(\mE^{\e}\left|\rho^{\e}_t(1)-\tilde{\rho}^{\e}_t(1)\right|^{2p}\Gamma_T^{\e}\right)^{1/2}\left(\mE^{\e}\left|\tilde{\rho}^{\e}_t(1)\right|^{-2p}\Gamma_T^{\e}\right)^{1/2}\\
&\leq&2^{p-1}\left(\mE^{\e}\left|\rho^{\e}_t(\phi)-\tilde{\rho}^{\e}_t(\phi)\right|^{4p}\right)^{1/4}\left(\mE^{\e}\left|\tilde{\rho}^{\e}_t(1)\right|^{-4p}\right)^{1/4}\left(\mE^{\e}\left|\Gamma_T^{\e}\right|^2\right)^{1/2}\\
&&+2^{p-1}\|\phi\|^{p}\left(\mE^{\e}\left|\rho^{\e}_t(1)-\tilde{\rho}^{\e}_t(1)\right|^{4p}\right)^{1/4}\left(\mE^{\e}\left|\tilde{\rho}^{\e}_t(1)\right|^{-4p}\right)^{1/4}\left(\mE^{\e}\left|\Gamma_T^{\e}\right|^2\right)^{1/2}\\
&\leq&\|\phi\|^{p}\frac{C}{\left[1-M\right]^p}\(\mE^\e\((2+2M)\|z_0\|_{\mH}+2 R(\omega)\)^{8p}\)^{1/8}\\
&&\(e^{-\frac{4\mu}{\e} tp}+\frac{\e}{4\mu p}\)^{1/4}\left(\mE^{\e}\left|\Gamma_T^{\e}\right|^2\right)^{1/2}\\
&=&\|\phi\|^{p}\frac{C}{\left[1-M\right]^p}\(\mE\((2+2M)\|z_0\|_{\mH}+2 R(\omega)\)^{8p}\left|\Gamma_T^{\e}\right|^{-1}\)^{1/8}\\
&&\(e^{-\frac{4\mu}{\e} tp}+\frac{\e}{4\mu p}\)^{1/4}\left(\mE^{\e}\left|\Gamma_T^{\e}\right|^2\right)^{1/2}\\
&\leq&\|\phi\|^{p}\frac{C}{\left[1-M\right]^p}\(\mE\((2+2M)\|z_0\|_{\mH}+2 R(\omega)\)^{16p}\)^{1/16}\\
&&\cdot\left(\mE|\Gamma_T^\e|^{-2}\right)^{1/16}(e^{-\frac{4\mu}{\e} t p}+\frac{\e}{4\mu p})^{1/4}\left(\mE^{\e}\left|\Gamma_T^{\e}\right|^2\right)^{1/2}.
\de
In the following, we estimate $\mE|\Gamma_T^\e|^{-2}, \mE^{\e}\left|\Gamma_T^{\e}\right|^2$. By simple calculations, it holds that
\ce
\mE(\Gamma_T^\e)^{-2}&=&\mE\left(\exp\left\{-2\int_0^T \<h(x_s^\e, y_s^\e),\dif
W_3(s)\>_{\tilde{\mathbb{H}}^3}+\frac{2}{2}\int_0^T\|h(x_s^\e, y_s^\e)\|_{\mH^3}^2\dif s\right\}\right)\\
&=&\mE\Bigg[\left(\exp\left\{-2\int_0^T \<h(x_s^\e, y_s^\e),\dif W_3(s)\>_{\tilde{\mathbb{H}}^3}-\frac{2^2}{2}\int_0^T\|h(x_s^\e,
y_s^\e)\|_{\mH^3}^2\dif s\right\}\right)\\
&&\cdot \exp\left\{\frac{2^2}{2}\int_0^T\|h(x_s^\e, y_s^\e)\|_{\mH^3}^2\dif s+\frac{2}{2}\int_0^T\|h(x_s^\e,
y_s^\e)\|_{\mH^3}^2\dif s\right\}\Bigg]\\
&\leq&\exp\left\{3C_h^2T\right\},
\de
where the last step is based on the fact that $\exp\left\{-2\int_0^t \<h(x_s^\e, y_s^\e),\dif W_3(s)\>_{\tilde{\mathbb{H}}^3}-\frac{2^2}{2}\int_0^t\|h(x_s^\e,
y_s^\e)\|_{\mH^3}^2\dif s\right\}$ is an exponential martingale under
$\mP$.
By the similar deduction to above, we know that 
$$
\mE^{\e}\left|\Gamma_T^{\e}\right|^2\leq\exp\left\{C_h^2T\right\}.
$$
Thus,
\ce
\mE|\pi^{\e}_t(\phi)- \tilde{\pi}_t^\e(\phi)|^p&\leq&\|\phi\|^{p}\frac{C}{\left[1-M\right]^p}\(\mE\((2+2M)\|z_0\|_{\mH}+2 R(\omega)\)^{16p}\)^{1/16}\\
&&\cdot(e^{-\frac{4\mu}{\e} t p}+\frac{\e}{4\mu p})^{1/4}.
\de

Next, notice that there exists a countable algebra $\{\phi_i, i=1, 2, \cdots\}$ of $\cC_b^1(\mH^1)$ that strongly seperates points in $\mH^1$. By \cite[Theorem 3.4.5]{ek}, it furthermore holds that $\{\phi_i, i=1, 2, \cdots\}$ is convergence determining for the topology of weak convergence of probability measures. For two probability measures $\mu, \tau$ on $\mathscr{B}(\mH^1)$,  set
\ce
d(\mu, \tau):=\sum\limits_{i=1}^{\infty}\frac{|\int_{\mH^1}\phi_i\,\dif\mu-\int_{\mH^1}\phi_i\,\dif\tau|}{2^i},
\de
and then $d$ is a distance in the space of probability measures on $\mathscr{B}(\mH^1)$. Since $\{\phi_i, i=1, 2, \cdots\}$ is convergence determining for the topology of weak convergence of probability measures, $d$ induces the weak convergence. The proof is completed.
\end{proof}

\br
By the above theorem, we know that when $\e$ goes to zero, $\tilde{\pi}_t^\e$ approximates $\pi^{\e}_t$. Therefore, $\tilde{\pi}_t^\e$ could be understood as the nonlinear filtering problem for $\tilde{x}_t^\e$ with respect to $\mathcal{R}_t^\e$.
\er

\section{An example}\label{exam}

\bx
Let $D$ be a domain in $\mR^3$ with smooth boundary $\partial D$. Consider the following coupled hyperbolic and parabolic equation 
\be
&&v_{tt}+\gamma v_t-\Delta v=f(v,v_t, \theta)+\sigma_1\dot{\cW}_1, t>0, x\in D,\no\\
&&\theta_t-\frac{1}{\e}\kappa\Delta\theta=\frac{1}{\e}g(v,v_t, \theta)+\frac{\sigma_2}{\sqrt{\e}}\dot{\cW}_2, t>0, x\in D,\no\\
&&v=0, \quad \theta=0, \qquad t>0, \quad x\in\partial D,
\label{exeq}
\ee
where $\gamma\geq 0, \kappa>0$ are constants, $\Delta$ is the Laplace operator and $f: \mR^3\rightarrow\mR$ and $g: \mR^3\rightarrow\mR$ are Lipschitz continuous with a Lipschitz constant $L>0$.

The type of equations are usually used to describe a thermoelastic phenomenon in a random medium(c.f.\cite{ccl}). Here $v$ denotes the displacement and $\theta$ is the temperature. And the parameter $\gamma$ describes resistance forces, and the white noise processes $\cW_1$ and $\cW_2$ model random fluctuations in external loads ($\dot{\cW}_1$) and in thermal sources ($\dot{\cW}_2$). If the temperature evolves fastly, then the hyperbolic equation is coupled to a parabolic equation with different characteristic timescales.

Next, we rewrite Eq.(\ref{exeq}) as
\ce\left\{\begin{array}{l}
\dot{x}^\e=A{x}^\e+F(x^\e,y^\e)+\sigma_1 \dot{W_1}, \\
\dot{y}^\e=\frac{1}{\e}B{y}^\e+\frac{1}{\e}G(x^\e, y^\e)+\frac{\sigma_2}{\sqrt{\e}}\dot{W_2},
\end{array}
\right.
\de
where
\ce 
&&x^\e=\left(\begin{array}{cccc}
v\\
v_t
\end{array}\right), A=\left(\begin{array}{cccc}
0 & 1\\
\Delta & -\gamma
\end{array}\right), F(x^\e, y^\e)=\left(\begin{array}{cccc}
0\\
f(v,v_t, \theta)
\end{array}\right), W_1=\left(\begin{array}{cccc}
0\\
\cW_1
\end{array}\right),\\
&&y^\e=\theta, B=\kappa\Delta, G(x^\e, y^\e)=g(v,v_t, \theta), W_2=\cW_2,\\
&&
\de
We take $\mH^1=H_0^1(D)\times L_2(D)$ and $\mH^2=L_2(D)$, where $L_2(D)$ and $H_0^1(D)$ are the usual Sobolev spaces. Thus, Eq.(\ref{exeq}) is in our framework. Moreover, by simple calculation, we know that $(\bf{H_1})$--$(\bf{H_3})$ are satisfied with $\gamma_1=\gamma, \gamma_2=\kappa$. If $\kappa>L$, it follows from Theorem \ref{solexiuni} that for $z_0=(x_0, y_0)\in\mH:=\mH^1\times\mH^2$, Eq.(\ref{exeq}) has a unique mild solution $z^\e(t, 0, \omega; z_0)$ for $t\in[0,T]$ and $\omega\in\Omega$. And Theorem \ref{Reduction-Theo} admits us to obtain that for $0<\e\leq\e_0$, there exists the following reduced system on the random invariant manifold $\cM^\e$:
\ce\left\{\begin{array}{l}
\dot{\tilde{x}}^\e=A\tilde{x}^\e+F\left(\tilde{x}^\e,\tilde{y}^\e\right)+\sigma_1 \dot{W_1},\\
\tilde{y}^\e=H^\e( \theta_{\cdot}\omega,
\tilde{x}^\e).
\end{array}
\right.
\de

In the following, we consider the nonlinear filtering problem of Eq.(\ref{exeq}). Taking $\mH^3=L_2(D)$,  one can construct a cylindrical Brownian motion $W_3$ on $\mH^3$. And then we give an observation system by
\ce
r^{\e}_t= W_3(t)+\int_0^t \sin (x^\e_s) \dif s, \quad t\in[0,T].
\de
Thus, it is easy to justify that $h(x,y)=\sin x$ satisfies $(\bf{H_5})$. By Theorem \ref{filcon}, the distance of the nonlinear filtering for $x^\e$ and ``the nonlinear filtering" for $\tilde{x}^\e$ is characterized.
\ex

\section{Conclusions}\label{Con}

In the paper, we consider an effective filtering for a type of slow-fast data assimilation systems in general Hilbert spaces. It is proved that the origin system can be reduced to a system on a random invariant manifold constructed by the direct way. And then, the nonlinear filtering of the origin system can be approximated by that of the reduction system. 

Of course, one could consider other different expressions for the noisy terms in the equations (multiplicative, cylindrical, etc). We will investigate the possibility of doing a similar reduction and the related nonlinear filtering in the future.

\section{Appendix}\label{App}

In this Appendix, we present the proof of Theorem \ref{solexiuni}.

{\bf Proof of Theorem \ref{solexiuni}.} First of all, we prove that the system (\ref{slfasy}) has a unique mild solution. Define an operator $\cJ: \cC([s,s+T], \mH)\rightarrow\cC([s,s+T], \mH)$ by
\ce
\cJ(z^\e)(t):=\left(\begin{array}{cccc}
\cJ_1(z^\e)(t)\\
\cJ_2(z^\e)(t)
\end{array}\right):=\left(\begin{array}{cccc}e^{A(t-s)}x_0+\int_s^te^{A(t-r)}F(x_r^\e,y_r^\e)dr+\int_s^te^{A(t-r)}\sigma_1\dif W_1(r)\\
e^{\frac{B}{\e}(t-s)}y_0+\int_s^te^{\frac{B}{\e}(t-r)}\frac{1}{\e}G(x_r^\e,y_r^\e)dr+\int_s^te^{\frac{B}{\e}(t-r)}\frac{\sigma_2}{\sqrt{\e}}\dif W_2(r)\end{array}\right).
\de
And then $\cJ$ is well defined. In fact, based on ($\bf{H_1}$) and \cite[Theorem 5.2]{dpz}, $e^{A(t-s)}x_0$ and $\int_s^te^{A(t-r)}\sigma_1\dif W_1(r)$ are strongly continuous. Taking $t_1, t_2\in[s,s+T], t_1<t_2$, we obtain that 
\ce
&&\left\|\int_s^{t_1}e^{A(t_1-r)}F(x_r^\e,y_r^\e)dr-\int_s^{t_2}e^{A(t_2-r)}F(x_r^\e,y_r^\e)dr\right\|_{\mH^1}\\
&\leq&\left\|\int_s^{t_1}e^{A(t_1-r)}F(x_r^\e,y_r^\e)dr-\int_s^{t_1}e^{A(t_2-r)}F(x_r^\e,y_r^\e)dr\right\|_{\mH^1}\\
&&+\left\|\int_s^{t_1}e^{A(t_2-r)}F(x_r^\e,y_r^\e)dr-\int_s^{t_2}e^{A(t_2-r)}F(x_r^\e,y_r^\e)dr\right\|_{\mH^1}\\
&\leq&\int_s^{t_1}\left\|F(x_r^\e,y_r^\e)-e^{A(t_2-t_1)}F(x_r^\e,y_r^\e)\right\|_{\mH^1}dr+\int_{t_1}^{t_2}\|e^{A(t_2-r)}\|\|F(x_r^\e,y_r^\e)\|_{\mH^1}dr\\
&\leq&\int_s^{t_1}\left\|F(x_r^\e,y_r^\e)-e^{A(t_2-t_1)}F(x_r^\e,y_r^\e)\right\|_{\mH^1}dr+L\(\sup\limits_{t\in[s,s+T]}\|z_t^\e\|_{\mH}\)(t_2-t_1),
\de
where the last step is based on (\ref{lingro}) and (\ref{opecon2}). The dominated convergence theorem admits us to get that $\int_s^te^{A(t-r)}F(x_r^\e,y_r^\e)dr$ is strongly continuous. Thus, $\cJ_1(z^\e)(t)$ is strongly continuous. By the same deduction to above, we know that $\cJ_2(z^\e)(t)$ is also strongly continuous. So, $\cJ(z^\e)(t)$ is strongly continuous.

Next, we prove that $\cJ$ is contractive. For $z^{\e,1}, z^{\e,2}\in\cC([s,s+T], \mH)$, one can compute by ($\bf{H_1}$)-($\bf{H_3}$)
\ce
\sup\limits_{t\in[s,s+T]}\left\|\cJ_1(z^{\e,1})(t)-\cJ_1(z^{\e,2})(t)\right\|_{\mH^1}&=&\sup\limits_{t\in[s,s+T]}\left\|\int_s^te^{A(t-r)}\(F(x_r^{\e,1},y_r^{\e,1})-F(x_r^{\e,2},y_r^{\e,2})\)dr\right\|_{\mH^1}\\
&\leq&\sup\limits_{t\in[s,s+T]}\int_s^t\|e^{A(t-r)}\|\left\|F(x_r^{\e,1},y_r^{\e,1})-F(x_r^{\e,2},y_r^{\e,2})\right\|_{\mH^1}dr\\
&\leq&\sup\limits_{t\in[s,s+T]}\int_s^tL(\|x_r^{\e,1}-x_r^{\e,2}\|_{\mH^1}+\|y_r^{\e,1}-y_r^{\e,2}\|_{\mH^2})dr\\
&\leq&LT\sup\limits_{t\in[s,s+T]}\|z^{\e,1}_t-z^{\e,2}_t\|_{\mH},
\de
and
\ce
\sup\limits_{t\in[s,s+T]}\left\|\cJ_2(z^{\e,1})(t)-\cJ_2(z^{\e,2})(t)\right\|_{\mH^1}&=&\sup\limits_{t\in[s,s+T]}\left\|\int_s^te^{\frac{B}{\e}(t-r)}\frac{1}{\e}\(G(x_r^{\e,1},y_r^{\e,1})-G(x_r^{\e,2},y_r^{\e,2})\)dr\right\|_{\mH^2}\\
&\leq&\frac{1}{\e}\sup\limits_{t\in[s,s+T]}\int_s^t\|e^{\frac{B}{\e}(t-r)}\|\left\|G(x_r^{\e,1},y_r^{\e,1})-G(x_r^{\e,2},y_r^{\e,2})\right\|_{\mH^2}dr\\
&\leq&\frac{1}{\e}\sup\limits_{t\in[s,s+T]}\int_s^te^{-\frac{\gamma_2}{\e}(t-r)}L(\|x_r^{\e,1}-x_r^{\e,2}\|_{\mH^1}+\|y_r^{\e,1}-y_r^{\e,2}\|_{\mH^2})dr\\
&\leq&\frac{L}{\gamma_2}[1-e^{-\frac{\gamma_2}{\e}T}]\sup\limits_{t\in[s,s+T]}\|z^{\e,1}_t-z^{\e,2}_t\|_{\mH}.
\de
Thus,
\ce
\sup\limits_{t\in[s,s+T]}\left\|\cJ(z^{\e,1})(t)-\cJ(z^{\e,2})(t)\right\|_{\mH}&\leq&\sup\limits_{t\in[s,s+T]}\left\|\cJ_1(z^{\e,1})(t)-\cJ_1(z^{\e,2})(t)\right\|_{\mH^1}\\
&&+\sup\limits_{t\in[s,s+T]}\left\|\cJ_2(z^{\e,1})(t)-\cJ_2(z^{\e,2})(t)\right\|_{\mH^2}\\
&\leq&\(LT+\frac{L}{\gamma_2}\)\sup\limits_{t\in[s,s+T]}\|z^{\e,1}_t-z^{\e,2}_t\|_{\mH}.
\de
Taking $T_0$ such that $LT_0+\frac{L}{\gamma_2}<1$, we know that $\cJ$ is contractive. So, the system (\ref{slfasy}) has a unique mild solution $z^\e(t, s, \omega; z_0)$ for $t\in[s,s+T_0]$. If $T\leq T_0$, the proof is over; if $T>T_0$, one can easily extend the solution to the finite interval $[s,s+T]$ by considering $[s,s+T_0]$, $[s+T_0,s+2T_0]$, $[s+2T_0,s+3T_0]$ and so on.

Set 
$$
\varphi^\e(t,\omega)z_0:=z^\e(t, 0, \omega; z_0), \quad t\in\mR,
$$
and then it is easy to justify that for $s, t\in\mR$,
\ce
&&\varphi^\e(t+s,\omega)z_0=z^\e(t+s, 0, \omega; z_0)=z^\e(t+s, s, \omega; z^\e(s, 0, \omega; z_0)),\\
&&\varphi^\e(t,\t_s\omega)\varphi^\e(s,\omega)z_0=z^\e(t, 0, \t_s\omega; z^\e(s, 0, \omega; z_0)).
\de
Note that $W_1(r, \t^1_s\omega)=W_1(s+r, \omega)-W_1(s, \omega), W_2(r, \t^2_s\omega)=W_2(s+r, \omega)-W_2(s, \omega)$ for $r\in\mR, s\in\mR, \omega\in\Omega$ and $W_1(r, \t^1_s\cdot), W_2(r, \t^2_s\cdot)$ are still two-sided $\mH^1, \mH^2$-valued Brownian motions with covariance operators $K_1, K_2$, respectively. Thus, by uniqueness of the mild solution for the system (\ref{slfasy}), we know that $z^\e(t+s, s, \omega; z^\e(s, 0, \omega; z_0))=z^\e(t, 0, \t_s\omega; z^\e(s, 0, \omega; z_0))$ and $\varphi^\e(t+s,\omega)z_0=\varphi^\e(t,\t_s\omega)\varphi^\e(s,\omega)z_0$. That is, $\varphi^\e(t,\omega)$ is a RDS. The proof is completed.

\bigskip

\textbf{Acknowledgements:}

The author would like to thank Professor Xicheng Zhang for his valuable discussions.


\begin{thebibliography}{999}
\bibitem{la}
L. Arnold, \emph{Random Dynamical Systems}, Springer, Berlin, 1998.

\bibitem{BC}
A. Bain and D. Crisan, \emph{Fundamentals of Stochastic Filtering}, Springer, Berlin, 2009.

\bibitem{ccl} T. Caraballo, I. Chueshov, and J.A. Langa, Existence of invariant manifolds for coupled
parabolic and hyperbolic stochastic partial differential equations, {\it Nonlin.}, 18(2005)747-767.

\bibitem{ccc} T. Cass, M. Clarke and D. Crisan, The filtering equations revisited, In: Crisan D., Hambly B., Zariphopoulou T. (eds) \emph{Stochastic Analysis and Applications (2014)}. Springer Proceedings in Mathematics \& Statistics, vol 100.

\bibitem{cs} I. Chueshov and M. Scheutzow, Inertial manifolds and forms for stochastically perturbed retarded semilinear
parabolic equations, {\it J. Dyn. Diff. Eqns,} 13(2001)355-80.

\bibitem{dpz} G. Da Prato and J. Zabczyk: Stochastic Equations in Infinite Dimensions. Encyclopedia of Mathematics
and its Applications. Cambriddge: Cambridge University Press, 1992.

\bibitem{dls} J. Duan, K. Lu and  B. Schmalfuss, Invariant manifolds for stochastic partial differential equations, {\it Ann. Probab.}, 31(2003), 2109-2135.

\bibitem{ek} S. N. Ethier and T. G. Kurtz, {\em Markov Processes: Characterization and Convergence}. John Wiley \& Sons, 1986.

\bibitem{Fu}  H. Fu, X. Liu and J. Duan,
Slow manifolds for multi-time-scale stochastic evolutionary systems,
\emph{Comm. Math. Sci.},
11(2013), 141-162.

\bibitem{hs} S. Hobbs and S.S. Sritharan,  Nonlinear filtering theory for stochastic reaction-diffusion equations. In {\it Probability and Modern Analysis}. Gritsky, N., Goldstein, J., and Uhl, J.J. eds., Marcel Dekker, New York, 1996, pp. 219-234. 

\bibitem{Imkeller} P. Imkeller, N. S. Namachchivaya, N. Perkowski and H. C. Yeong,
 Dimensional reduction in nonlinear filtering: a homogenization approach,
\emph{The Annals of Applied Probability},
23(2013), 2290-2326.

\bibitem{RG}
R. Z. Khasminskii and G.Yin,
On transition densities of singularly perturbed diffusions with fast and slow components,
\emph{SIAM J. Appl. Math.},
56(1996), 1794-1819.

\bibitem{ps1} A. Papanicolaou and K. Spiliopoulos: Dimension reduction in statistical estimation of partially observed multiscale processes, {\it SIAM J. on Uncertainty Quantification,} 5(2017)1220-1247.

\bibitem{Park1}  J. H. Park, N. S. Namachchivaya and H. C. Yeong,
Particle filters in a multiscale environment: Homogenized hybrid particle filter,
\emph{J. Appl. Mech.},
78(2011), 1-10.

\bibitem{Park2}  J. H. Park, R. B. Sowers and N. S. Namachchivaya,
Dimensional reductionin nonlinear filtering,
\emph{ Nonlinearity},
23(2010), 305-324.

\bibitem{Park3}
J. H. Park, B. Rozovskii and R. B. Sowers,
Efficient nonlinear filtering of a singularly perturbed stochastic hybrid system,
\emph{LMS Journal of Computation and Mathematics},
14(2011), 254-270.

\bibitem{andrew} G. A. Pavliotis and A. M. Stuart, Multiscale methods: averaging and homogenization, Springer Science+Business Media, New York, 2008.

\bibitem{q1} H. Qiao, Convergence of nonlinear filtering for stochastic dynamical systems with L\'evy noises, https://arxiv.org/abs/1707.07824.

\bibitem{q2} H. Qiao, Y. Zhang and J. Duan, Effective filtering on a random slow manifold, {\it Nonlinearity,} 31(2018) 4649-4666.

\bibitem{RBL}
B. L. Rozovskii,
\emph{Stochastic Evolution System: Linear Theory and Application to nonlinear Filtering},
Springer, New York, 1990.

\bibitem{ZS}
Z.  Schuss,
\emph{Nonlinear Filtering and Optimal Phase Tracking},
Springer, New York, 2012.

\bibitem{Sch}
 B. Schmalfu{\ss} and  R. Schneider,
 Invariant manifolds for random dynamical systems with slow and fast variables,
\emph{J. Dyna. Diff. Equa.},
20(2008), 133-164.

\bibitem{sr} S. S. Sritharan, Nonlinear filtering of stochastic Navier-Stokes equations. In {\it Nonlinear Stochastic PDEs: Burgers Turbulance and Hydrodynamic Limit}. Funaki, T., and Woyczynski, W.A. eds., Springer-Verlag, New York, 1995, pp. 247-260.

\bibitem{yz} Y. Ying and F. Zhang, Potentials in Improving Predictability of Multiscale Tropical Weather Systems Evaluated through Ensemble Assimilation of Simulated Satellite-Based Observations, {\it Journal of the Atmospheric Sciences}, 2018, https://doi.org/10.1175/JAS-D-17-0245.1.


\bibitem{zqd} Y. Zhang, H. Qiao and J. Duan, Effective filtering analysis for non-Gaussian dynamic systems, appear in {\it Applied Mathematics and Optimization}.
\end{thebibliography}
\end{document}